\documentclass[review,12pt]{elsarticle}

\usepackage{setspace}
\usepackage{subcaption}
\usepackage{multirow}
\usepackage{multicol}
\usepackage{graphicx}
\RequirePackage{amsmath} 
\usepackage{bm}
\usepackage{xfrac}
\usepackage[most]{tcolorbox}
\usepackage{float}

\usepackage{amsmath,esint}

\usepackage{lineno,hyperref}

\journal{Applied Mathematical Modelling}

\biboptions{sort&compress}

\begin{document}

\begin{frontmatter}

\title{Finite Volume Simulation Framework for Die Casting with Uncertainty Quantification}

\author{Shantanu Shahane*\fnref{Corresponding Author}}
\author{Narayana Aluru} \author{Placid Ferreira} \author{Shiv G Kapoor} \author{Surya Pratap Vanka}
\address{Department of Mechanical Science and Engineering\\
	University of Illinois at Urbana-Champaign \\
	Urbana, Illinois 61801}
\fntext[Corresponding Author]{Corresponding Author}

%

\begin{abstract}
The present paper describes the development of a novel and comprehensive computational framework to simulate solidification problems in materials processing, specifically casting processes. Heat transfer, solidification and fluid flow due to natural convection are modeled. Empirical relations are used to estimate the microstructure parameters and mechanical properties. The fractional step algorithm is modified to deal with the numerical aspects of solidification by suitably altering the coefficients in the discretized equation to simulate selectively only in the liquid and mushy zones. This brings significant computational speed up as the simulation proceeds. Complex domains are represented by unstructured hexahedral elements. The algebraic multigrid method, blended with a Krylov subspace solver is used to accelerate convergence. State of the art uncertainty quantification technique is included in the framework to incorporate the effects of stochastic variations in the input parameters. Rigorous validation is presented using published experimental results of a solidification problem.
\end{abstract}


\end{frontmatter}
\section{Introduction}
Among the various processes used by the industry to produce a finished part, die casting is an important manufacturing process in which liquid metal is injected into a die and solidified. Usually die casting is used for aluminum and magnesium alloys with steel molds. Die cast products are more commonly used in automotive and housing industries. Several complex processes involving a large number of process parameters affect the final product quality. Due to the recent advances in computing hardware and software, it is now possible to simulate the physics of these processes using numerical simulations. These simulations provide detailed flow and temperature histories and can be used to estimate the final product quality. Frequently, it is difficult to measure and tightly control the process parameters like initial melt temperature, mold temperature and alloy component concentration. However, they can have a significant impact on the process as well as the predicted product strength.
\par Temperature evolution and velocity distribution during solidification of pure metals and metal alloys has been analyzed by many researchers. \citet{bennon1987continuum_1} used a continuum mixture model to solve the momentum and energy equations and applied it to the solidification problem of a rectangular cavity filled with a binary aqueous solution \cite{bennon1987continuum_2}. From their work, it is clear that although the continuum formulation has a limitation of smearing the interface, it is efficient than the multiple region method in which the governing equations are solved in each phase separately with appropriate interface conditions. Implementing multiple region method for a practical complex geometry with irregular interfaces is quite difficult and computationally expensive. \citet{voller1987fixed} used the Darcy's law of drag and the enthalpy method to model the mushy zone during solidification of a square cavity. Two vertical walls were held at below and above the melting point and the horizontal walls were thermally insulated. The damping of velocities in the mushy zone due to Darcy's law is clearly seen in the velocity vector plot. \citet{vynnycky2007analytical} solved a similar problem on a rectangular enclosure analytically. They used an asymptotic approach to solve the non-dimensional momentum and energy equations for the case of the Rayleigh and Stefan numbers much larger and smaller than unity, respectively. Comparison with finite element transient numerical simulation showed that such an asymptotic simplification is possible. \citet{bennett2007adaptive} used adaptive grid with local rectangular refinement (LRR) to model solidification with fluid flow in multiple rectangular geometries. The adaptive LRR grid needed half the storage and computational effort compared to the non-adaptive grid. \citet{wang2010comprehensive} discussed a numerical model for melting in a rectangular cavity with natural convection at high Rayleigh number ($10^8$). They used the consistent update technique (CUT) algorithm with a multigrid solver. \citet{plotkowski2015estimation} simulated analytically and numerically heat transfer and fluid flow for solidification in a rectangular cavity. They used scaling analysis to simplify the governing equations of the mixture model followed by an analytical solution and comparison with a finite volume solution for an Al-Cu binary alloy. \citet{hu2017lattice} studied the three-dimensional phase change problem with natural convection using the Lattice Boltzmann method. They simulated melting in a cubical cavity and in cavities with inner rectangular cylinders and sphere. \citet{cleary2010short} used smoothed particle hydrodynamics to model filling and solidification in high pressure die casting. They simulated practical industrial case studies like differential cover, an electronic housing and a door lock plate. \citet{gau1986melting} performed experiments to understand the importance of natural convection during solidification and melting of pure metal. \citet{quillet2007benchmark} and \citet{hachani2012experimental} reported the temperature distribution and macro-segregation of solute with an experimental solidification of a Tin based alloy in an ingot. Such experimental work is useful for calibration and validation of the numerical methods.
\par Solidification phenomena of die casting involves interplay between heat transfer and flow due to natural convection. All the references discussed above simulated solidification on rectangular or cuboidal geometries and hence, structured Cartesian grids were used. Practical die casting geometries are complex and thus, unstructured hexahedral elements are utilized in this work. Multiple algorithms are discussed in the literature to simulate fluid flow on unstructured grids \cite{davidson1996pressure,haselbacher1999grid,kim2000second,ferziger2012computational,mathur1997pressure,muzaferija1997finite}. Strategy implemented in this research is an extension of the work by \citet{mathur1997pressure} and \citet{muzaferija1997finite}. Solution of the fluid flow and the pressure Poisson equation involves a significant computational effort. As the solidification problem proceeds, the volume of liquid zone deceases continuously and it is not necessary to solve for fluid flow in the solid zone. We therefore, have developed a consistent procedure to modify the coefficients of the discretized momentum and pressure Poisson equations in order to satisfy the mass continuity equation. Algebraic multigrid \cite{yang2002boomeramg,ruge1987algebraic} and Krylov subspace solvers from the open source library HYPRE \cite{HYPRE_website} are used to solve the linear systems of equations.
\par For die casting, the microstructure parameters like grain size and dendritic arm spacing are important as they affect the final product quality. Phase field modeling \cite{provatas1999adaptive,karma1998quantitative,ramirez2004phase} is a popular method used to study the evolution of the microstructure during solidification. The phase field method simulates the growth of each dendrite and thus, it is computationally expensive at the length scale of die cast products. In this research, an empirical relation from the work of \citet{backer2007microporosity} is used to estimate secondary dendritic arm spacing. There are various models \cite{maxwell1975simple,greer2000modelling,desnain1990prediction} suggested for grain size estimation during solidification. Here, the isothermal crystal growth model \cite{greer2000modelling} is used. 
\par Use of deterministic simulations alone to analyze the engineering systems is incomplete due to the lack of precisely defined input data. Thus, there has been a growing interest \cite{xiu2002modeling,xiu2003modeling,knio2001stochastic,carnevale2013uncertainty,marepalli2014quantifying,ganapathysubramanian2007sparse,fezi2017uncertainty,hosder2006non,kumar2016efficient,fajraoui2017analyzing} in coupling uncertainty propagation techniques with the deterministic numerical simulations to estimate the effects of stochastic variations in the input process parameters on the outputs. The polynomial chaos expansion is a popular method used to estimate the relation between input and output parameters. Stochastic Galerkin projection \cite{xiu2002modeling,xiu2003modeling,knio2001stochastic} and collocation \cite{carnevale2013uncertainty,marepalli2014quantifying,ganapathysubramanian2007sparse} are two strategies to estimate the coefficients of the polynomial chaos expansion. Stochastic Galerkin method is an intrusive method since it requires solution of a new set of equations and thus, modification of the underlying deterministic code which becomes a significant additional effort. Hence, recently non-intrusive stochastic collocation methods have been developed which need multiple evaluations of the deterministic simulation at predefined collocation points obtained by sampling from the probability distribution function of the input parameters. Values of outputs estimated at these samples are then used to estimate the coefficients of the polynomial chaos expansion.
\par In this paper, the traditional fractional step approach is modified to deal with the additional terms in the Navier-Stokes equation due to solidification. The discretized system of equations is altered in a way so as to simulate only in the liquid zone thus, increasing computational efficiency. In order to incorporate the effects of stochastic variations in the input process parameters, a parameter uncertainty propagation module has been developed in conjunction with the deterministic simulations. The method of polynomial chaos expansion is used to estimate the relation between input and output parameters and stochastic collocation is used as a wrapper on the underlying deterministic simulation. The framework has been validated against published experimental results followed by demonstration for solidification of two complex geometries.
\section{Numerical Model Description}
\subsection{Governing Equations}
Solidification, heat transfer and fluid flow due to natural convection are modeled. It is assumed that there is no macro-segregation during solidification and the metal is solidified at nominal composition. It is further assumed that the solutes do not contribute to buoyancy and the flow is incompressible. Thus, the set of governing equations consists of the standard Navier-Stokes equations with additional terms for solidification \cite{plotkowski2015estimation}. They are:
\begin{equation} 
	\nabla \cdot \textbf{u} =0
	\label{Eq:Continuity}
\end{equation}
\begin{equation} 
	\rho \frac{\partial \textbf{u}}{\partial t} + \nabla \cdot (\rho \textbf{u} \otimes \textbf{u}) = \nabla \cdot ( \mu \nabla \textbf{u}) - \nabla P - \frac{\mu}{K}\textbf{u} - \textbf{g} \rho \beta (T-T_{ref})
	\label{Eq:Momentum}
\end{equation}
\begin{equation} 
	K= \frac{\lambda^2(1-f_s)^3}{180 f_s^2}
	\label{Eq:Darcy_drag}
\end{equation}
\par Here, $\textbf{u}$ is the velocity vector, $\rho$ is density, $t$ is time, $\mu$ is dynamic viscosity, $\textbf{g}$ is gravity vector, $\beta$ is coefficient of thermal expansion, $P$ is pressure, $K$ is isotropic permeability of the dendritic array, $\lambda$ is dendrite arm spacing and $f_s$ is solid fraction.
\par To model the effects of natural convection, the Boussinesq approximation is used. This is a valid assumption for problems with moderate density variations in the domain. The fluid is modeled as a constant density fluid except for the additional buoyancy term $-\textbf{g }\rho \beta (T-T_{ref})$ in the momentum equation (\ref{Eq:Momentum}) \cite{spiegel1960boussinesq}.
\par The Darcy drag term ($ \frac{\mu}{K}\textbf{u} $) represents increased resistance to the flow in the mushy zone. We have used the Blake-Kozeny model (Eq.~(\ref{Eq:Darcy_drag})) which estimates the isotropic permeability ($K$) of the dendritic array. In the liquid region, solid fraction is zero and permeability tends to infinity making the Darcy drag term to go to zero. When solid fraction is unity, permeability tends to zero and thus the coefficient of Darcy drag term goes to infinity. For stability, this coefficient is added to the diagonal term of the discretized momentum equations. As a result, the velocities in the solid region go to zero. In the mushy zone, the drag term reduces the velocities compared to the liquid zone.
\par The energy equation is written in terms of temperature as:
\begin{equation} 
	\rho C_p \frac{\partial T}{\partial t} + \rho C_p (\nabla \cdot \textbf{u} T) = \nabla \cdot ( k \nabla T) + \rho L_f \frac{\partial f_s}{\partial t} 
	\label{Eq:Energy}
\end{equation}
where,
\begin{equation}
	f_s(T)=
	\begin{cases}
		0 & \text{if } T  >T_{liq}\\
		1 & \text{if } T  <T_{sol}\\
		1 - \left( \frac{T-T_f}{T_{liq} - T_f}\right)^{\frac{1}{k_p-1}}  & \text{otherwise}
	\end{cases}
	\label{Eq:Solid_fraction}
\end{equation}
\par Here, $T$ is temperature, $C_p$ is specific heat, $k$ is thermal conductivity, $L_f$ is latent heat of fusion, $k_p$ is partition coefficient, $T_f$ is freezing temperature and $T_{liq}$ is liquidus temperature. To close the system, the Gulliver-Scheil equation (\ref{Eq:Solid_fraction}) \cite{dantzig2001modeling} for the solid fraction is used.
\subsection{Solution Algorithm}\label{Sec:Fractional Step Algorithm}
We have developed a new software OpenCast in an object oriented C++ environment. A finite volume method on a collocated grid is used to discretize the governing equations. The fractional step method \cite{harlow1965numerical}, modified to account for the solidification is used to integrate the equations. First, the momentum equations (\ref{Eq:Momentum_frac_step_u*}) without the pressure gradient term is solved to estimate an intermediate velocity ($\textbf{u*}$) field. The buoyancy and Darcy drag terms are included in this step. Second order accurate Crank-Nicolson scheme for the diffusion terms and the Adams-Bashforth scheme for the convection terms are used for temporal discretization. The coefficient of the drag term ($\frac{\mu}{K}$) in the mushy zone is treated fully implicitly as:
\begin{equation} 
	\rho \frac{\textbf{u*} - \textbf{u}^n}{\Delta t} + \frac{\mu}{K}\textbf{u*} = -Conv(\textbf{u}^n,\textbf{u}^{n-1}) + Diff(\textbf{u*},\textbf{u}^n) + Buoy(T^n)
	\label{Eq:Momentum_frac_step_u*}
\end{equation}
where, the operators $Conv$, $Diff$ and $Buoy$ represent the discretized convection, diffusion and buoyancy terms respectively. The full momentum equation is similarly discretized with an implicit pressure gradient term, given as:
\begin{equation} 
	\rho \frac{\textbf{u}^{n+1} - \textbf{u}^n}{\Delta t} + \frac{\mu}{K}\textbf{u*} = -Conv(\textbf{u}^n,\textbf{u}^{n-1}) + Diff(\textbf{u*},\textbf{u}^n) - (\nabla P)^{n+1} + Buoy(T^n)
	\label{Eq:Momentum_frac_step_u}
\end{equation}
Subtracting equation (\ref{Eq:Momentum_frac_step_u*}) from (\ref{Eq:Momentum_frac_step_u}) gives the velocity correction equation.
\begin{equation} 
	\textbf{u}^{n+1}=\textbf{u*}- (\nabla P)^{n+1} \frac{\Delta t}{\rho}
	\label{Eq:Velocity_correction}
\end{equation}
Taking divergence of the velocity correction equation (\ref{Eq:Velocity_correction}) and invoking the continuity constraint gives the equation for pressure:
\begin{equation} 
	 \nabla \cdot \left(\frac{\nabla P }{\rho}\right)^{n+1} = \frac{\nabla \cdot \textbf{u*}}{\Delta t} 
	\label{Eq:Pressure_poisson}
\end{equation}
The overall solution algorithm to advance from time-step $n$ to $n+1$ is as follows:
\begin{enumerate}
	\item Solve for $\textbf{u*}$ using equation (\ref{Eq:Momentum_frac_step_u*}). Since the diffusion term is implicit, solution is obtained iteratively
	\item Solve the pressure Poisson equation (\ref{Eq:Pressure_poisson}) iteratively to estimate $P^{n+1}$
	\item Correct the velocities ($\textbf{u}^{n+1}$) using equation (\ref{Eq:Velocity_correction})
	\item Solve the energy equation (\ref{Eq:Energy}) together with solid fraction relation (Eq.~(\ref{Eq:Solid_fraction})) to estimate temperature and solid fraction
	\item Estimate micro-structure parameters such as grain size and yield strength using the empirical relations
\end{enumerate}
\subsection{Discretization on an Unstructured Grid}
\par Practical die casting geometries are quite complex. Cartesian grids introduce high stair-casing errors near the boundaries. Thus, OpenCast uses unstructured grids with tetrahedral and hexahedral finite volumes. First, a tetrahedral mesh is generated using the open source software GMSH \cite{geuzaine2009gmsh}. Tetrahedral elements have higher cross diffusion terms due to mesh skewness compared to hexahedral elements. Thus, we further use an open source tool TETHEX \cite{Tethex_github} to subdivide the tetrahedrons to hexahedrons.
\par The governing equations (\ref{Eq:Continuity}), (\ref{Eq:Momentum}) and (\ref{Eq:Energy}) can be written as a transport equation for a general scalar $\phi$:
\begin{equation} 
\rho \frac{\partial \phi}{\partial t} + \nabla \cdot(\rho \textbf{u} \phi) = \nabla \cdot ( \Gamma \nabla \phi) + S_{\phi}
\label{Eq:scalar_transport}
\end{equation}
where, $\phi$ is any scalar field, $\Gamma$ is the diffusion coefficient, and $S_{\phi}$ is the source term.
\begin{figure}[h]
	\centering
	\includegraphics[width=0.63\textwidth]{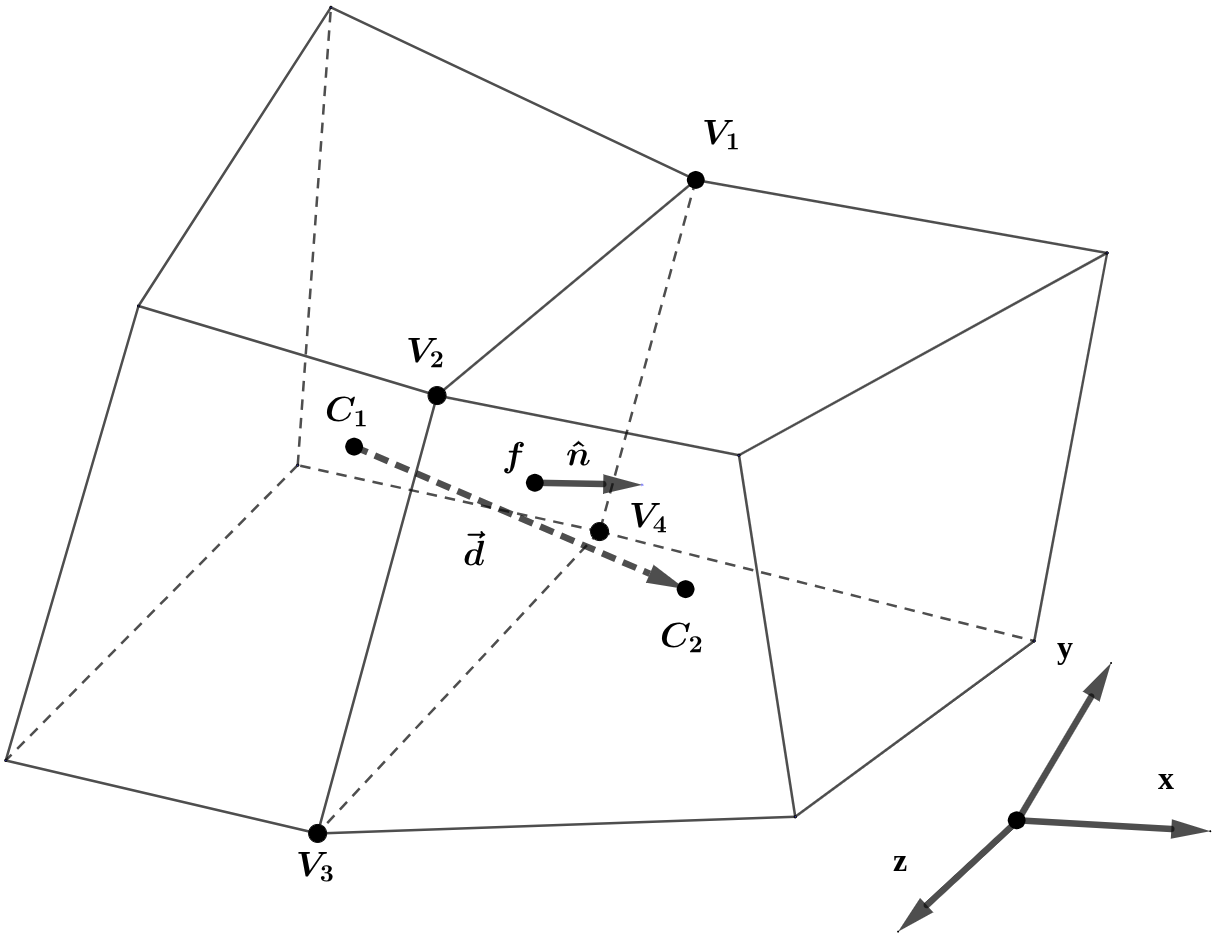}
	\caption{Unstructured Hexahedral Control Volumes}
	\label{Fig:Hexahedral_CV}
\end{figure}
\par Figure~\ref{Fig:Hexahedral_CV} shows two adjacent hexahedral control volumes sharing a common face with vertices $V_1$, $V_2$, $V_3$ and $V_4$. $C_1$ and $C_2$ are cell centers and $f$ is the face center. $\hat{n}$ is the unit vector normal to face and in an outward direction with respect to cell $C_1$. $\vec{d}$ is the distance vector from $C_1$ to $C_2$. We use a collocated finite volume formulation with all the field variables stored at cell centers.
\par The surface integral of the diffusion term is approximated as a summation over all the six faces of the cell. The inner product of the normal and the face centered gradient at each face is split into two terms \cite{mathur1997pressure}:
\begin{equation}
\begin{aligned}
\hat{n} \cdot \nabla \phi \Bigr|_{\substack{f}}=&\left(\frac{\vec{d} \cdot \nabla \phi \Bigr|_{\substack{f}}}{\hat{n} \cdot \vec{d}}\right) - \left(\frac{\vec{d} \cdot \nabla \phi \Bigr|_{\substack{f}}}{\hat{n} \cdot \vec{d}} - \hat{n} \cdot \nabla \phi \Bigr|_{\substack{f}}\right)\\
=&\left(\frac{\phi_{C_2} - \phi_{C_1}}{\hat{n} \cdot \vec{d}}\right) + \left(\hat{n} - \frac{\vec{d}}{\hat{n} \cdot \vec{d}} \right) \cdot \nabla \phi \Bigr|_{\substack{f}}
\end{aligned}
\label{Eq:Diffusion_2}
\end{equation}
The first and second terms of equation (\ref{Eq:Diffusion_2}) are direct and cross diffusion terms, respectively. For a structured grid, $\hat{n}$ is parallel to $\vec{d}$ and the cross diffusion term is identically zero as the direct diffusion term reduces to the central difference approximation of first derivative at face center.
\par In order to estimate the face centered gradient, the strategy used by \citet{mathur1997pressure} for two dimensional grids is extended. A local co-ordinate system is defined with $\xi: (C_1 C_2)$, $\eta: (V_1 V_3)$ and $\zeta: (V_2 V_4)$ (Fig.~\ref{Fig:Hexahedral_CV}) as the three axes. $x$, $y$ and $z$ are the axes of the global frame of reference. The gradients in both these frames are related by the chain rule of differentiation.
\begin{equation} 
\begin{bmatrix}
\phi_{\xi}\\
\phi_{\eta} \\
\phi_{\zeta}
\end{bmatrix} =
\begin{bmatrix} x_{\xi} & y_{\xi} & z_{\xi}  \\ 
x_{\eta} & y_{\eta} & z_{\eta}\\
x_{\zeta} & y_{\zeta} & z_{\zeta} \end{bmatrix}
\begin{bmatrix}
\phi_x\\
\phi_y \\
\phi_z
\end{bmatrix} 
\label{Eq:Diffusion_3}
\end{equation}
where, the subscripts denote derivatives. The Jacobian matrix entries come from the co-ordinates of the cell centers and the vertices. Value of $\phi$ at each vertex is estimated by averaging from the neighboring cells of the vertex. Thus, the face centered gradient $ \nabla \phi \Bigr|_{\substack{f}}=[\phi_x, \phi_y, \phi_z]^T$ is estimated by inverting the Jacobian matrix in equation (\ref{Eq:Diffusion_3}) and multiplying by $\begin{bmatrix}	\phi_{\xi}, \phi_{\eta}, \phi_{\zeta} \end{bmatrix}^T$.
\par The surface integral of the convection term is approximated as a summation over all the six faces of the cell. The face value of the field $\phi$ is estimated by interpolating from the two neighboring cells which share the face. The volume flux passing through the face ($\hat{n} \cdot \textbf{u} \Delta A$) satisfies the discrete continuity equation. The cross diffusion term has to be accounted for in the computation of the volume flux. The details are given in Section~\ref{Sec:Face Volume Flux Computation and Velocity Correction}.
\subsection{Special Modifications in Solidification Regions}
For solidification problems, some additional steps are needed in order to handle the extra terms such as the Darcy drag and the latent heat terms in the momentum and energy equations respectively. The velocities in the solid region should go to zero and in the mushy zone, velocities should be significantly lower than the fully liquid region. Simultaneously, the continuity equation has to be satisfied by the face velocities for each control volume. Thus, special care has to be taken in the solution process of the pressure Poisson equation and the velocity correction step.
\par
\begin{figure}[h]
	\centering
	\includegraphics[width=0.5\textwidth]{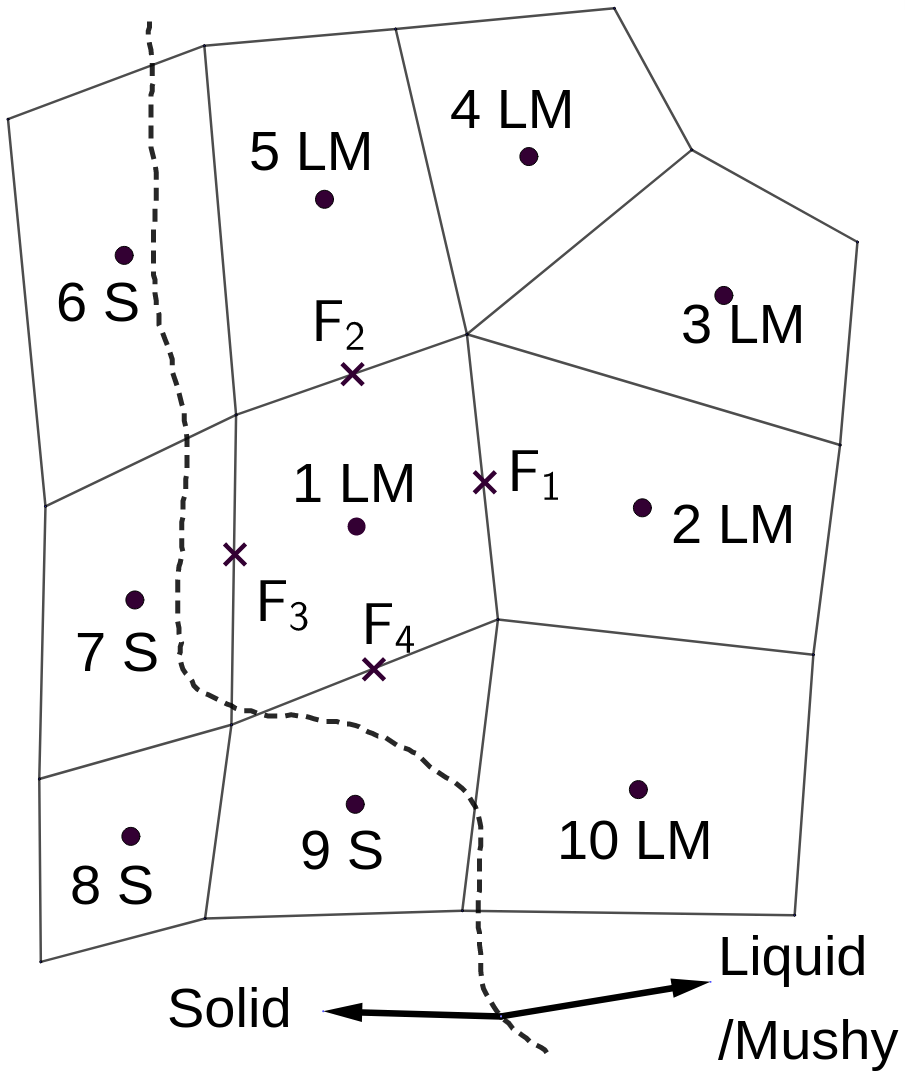}
	\caption{Control Volumes with Solid Liquid Interface}
	\label{Fig:Coeff_connectivity}
\end{figure}

 Figure~\ref{Fig:Coeff_connectivity} shows a typical distribution of phases during the solidification process. For ease of visualization, a two dimensional schematic is shown and the same idea has been generalized to three dimensions. The dotted line shows a solid-mushy zone interface. The control volumes (cells) on the left of the line are solid and on the right are either liquid or in mushy zone. Cells are labeled with tags S:Solid, LM: Liquid or Mushy. The S or LM tag is assigned to each cell based on the temperature at the previous time step and the liquidus and solidus temperatures of the alloy.
\subsubsection{Momentum equation}
As described in section \ref{Sec:Fractional Step Algorithm}, the modified momentum equations (\ref{Eq:Momentum_frac_step_u*}) are solved to estimate the intermediate velocities ($\textbf{u*}$). After discretization, the Darcy drag coefficient ($\frac{\mu}{K}$) is added to the diagonal term of the linear equations. In the pure liquid region, this coefficient is zero and thus, it does not have any effect. In the mushy zone, it is finite and non-zero and thus, it acts like a resistance to the flow. In the fully solidified region, it is a large number and thus, the $\textbf{u*}$ tends to zero. From computational efficiency point of view, it is not necessary to solve for $\textbf{u*}$ in the solidified cells since it is zero. Therefore, at each time step, before solving the linearized system of equations, the matrix rows corresponding to solidified cells are removed. As these solidified cells are connected to the neighboring mushy or liquid zone cells, the rows corresponding to the neighboring cells have to be modified in a consistent manner. Consider the row corresponding to cell number 1 in Fig.~\ref{Fig:Coeff_connectivity}: 
\begin{equation} 
\left[A_1, A_2, \dots, A_{10}\right] \left[\phi_1,\phi_2,\dots,\phi_{10}\right]^T = \left[S_1, S_2, \dots, S_{10}\right]^T
\label{Eq:Discrete_momentum_1}
\end{equation}
where, $\phi$ is any component of $\textbf{u*} = [u*,v*,w*]$. Originally, $\phi_1$ is connected to all the neighboring cells from 2 to 10. But since cells 6, 7, 8 and 9 are solidified and their velocity is zero, those rows and columns are deleted from equation (\ref{Eq:Discrete_momentum_1}). Thus, the reduced row becomes:
\begin{equation} 
\left[A_1, A_2, A_3, A_4, A_5, A_{10}\right] \left[\phi_1,\phi_2,\phi_3,\phi_4,\phi_5,\phi_{10}\right]^T = \left[S_1, S_2, S_3,S_4,S_5, S_{10}\right]^T
\label{Eq:Discrete_momentum_2}
\end{equation}
\subsubsection{Pressure Poisson Equation}\label{Sec:Pressure Poisson Equation}
Similar to the momentum equations, the matrix rows corresponding to solidified cells are deleted from the discrete pressure Poisson equation. But in this case, the rows of the neighboring liquid or mushy cells cannot be updated by just deleting the connections of the solid cells as Neumann boundary conditions have to be applied.
For any solidified cell, incoming or outgoing flow through all of its faces should be made zero. This is achieved as follows. Each face is shared by exactly 2 cells (face owners). All the cells which share a common vertex with a face are known as its neighbors. The face centered gradient is computed using the values at all connected neighboring cells. For example, cell numbers 1 and 2 are the owners of face F$_1$ whereas, cells 1, 2, 3, 4, 5, 9 and 10 are its neighbors. The following cases arise for each face:
\begin{enumerate}
	\item None of the neighbor cells is solidified: no change in the face centered gradient coefficient (eg. all the faces of cell number 3) is required
	\item None of the owners are solidified but at least one neighbor cell is solidified: flow through the face is allowed but the face centered gradient coefficient has to be modified as the solidified cells are removed from the linear set of equation (eg. faces F$_1$ and F$_2$)
	\item At least one owner is solidified: flow through the face is blocked; $\nabla P \cdot \hat{n} = 0$ thus, contribution of this face in the integrated diffusion term ($\oiint_S \Gamma \hat{n} \cdot \nabla P dS$) is zero (eg. faces F$_3$ and F$_4$)
\end{enumerate}
Consider the face F$_2$ for modification of face centered gradient coefficient. Original coefficients for gradient computation at face F$_2$ which are valid if none of its neighbor cells are solid are given by:
\begin{equation} 
\renewcommand\arraystretch{1.8}
\begin{bmatrix}
\frac{\partial P}{\partial x}\\
\frac{\partial P}{\partial y}
\end{bmatrix}_{F_2} \approx 
\begin{bmatrix} A_{x1} & A_{x2} & A_{x3} &A_{x4} & A_{x5} & A_{x6} & A_{x7}\\
 A_{y1} & A_{y2} & A_{y3} &A_{y4} & A_{y5} & A_{y6} & A_{y7} \end{bmatrix}
\begin{bmatrix}
P_1 & P_2 & P_3 & P_4 & P_5 & P_6 & P_7
\end{bmatrix} ^T
\label{Eq:Face_gradient}
\end{equation}
Since cell numbers 6 and 7 have solidified, their contribution has to be removed from equation (\ref{Eq:Face_gradient}). Thus, the last 2 columns are deleted and those coefficients are smeared equally in the remaining columns for ex., $A_{x1}$ is modified to $B_{x1} = A_{x1} + (A_{x6} + A_{x7})/5$ and $A_{y1}$ to $B_{y1} = A_{y1} + (A_{y6} + A_{y7})/5$. Since there are 5 cells remaining, the division by 5 is required. After modification, equation (\ref{Eq:Face_gradient}) becomes:
\begin{equation} 
\renewcommand\arraystretch{1.8}
\begin{bmatrix}
\frac{\partial P}{\partial x}\\
\frac{\partial P}{\partial y}
\end{bmatrix}_{F_2} \approx 
\begin{bmatrix} B_{x1} & B_{x2} & B_{x3} &B_{x4} & B_{x5}\\
B_{y1} & B_{y2} & B_{y3} &B_{y4} & B_{y5}\end{bmatrix}
\begin{bmatrix}
P_1 & P_2 & P_3 & P_4 & P_5
\end{bmatrix} ^T
\label{Eq:Face_gradient_modified}
\end{equation}
The steps for modification of the discretized pressure Poisson equation are as follows:
\begin{enumerate}
	\item Identify the faces with none of the owners solidified but at least one neighbor cell is solidified and smear the coefficients as described above
	\item If at least one owner of the face is solidified, set all of its coefficients to zero
	\item Loop over all the cells:
	\begin{itemize}
		\item If none of its face coefficients are modified, its coefficients do not change
		\item If at least one of its face coefficients is modified, re-assemble its coefficients
	\end{itemize}
\end{enumerate}
These steps remove the contribution of the solidified cells carefully and reduce the computational effort significantly.
\subsubsection{Face Volume Flux Computation and Velocity Correction}\label{Sec:Face Volume Flux Computation and Velocity Correction}
The collocated finite volume formulation uses the face centered volume fluxes ($\hat{n} \cdot \textbf{u} \Delta A$) in the continuity equation so as to avoid the checker-boarding of pressure. If there is an inconsistency in the numerical formulation of the pressure Poisson equation and the flux computations, there can be a gain or loss of mass and convergence problems. This section describes a consistent method used in the current code to handle solidification.
\par The volume flux is obtained by taking inner product of the velocity correction equation (\ref{Eq:Velocity_correction}) at the face center with face normal and multiplying by face area:
\begin{equation} 
	\hat{n} \cdot \textbf{u}^{n+1} \Delta A \Bigr|_{\substack{f}}= \hat{n} \cdot \textbf{u*} \Delta A \Bigr|_{\substack{f}} - \hat{n} \cdot (\nabla P)^{n+1} \Delta A \Bigr|_{\substack{f}} \frac{\Delta t}{\rho}
	\label{Eq:Volume_flux}
\end{equation}
$\textbf{u*}\Bigr|_{\substack{f}}$ is estimated by averaging the cell values from the two owner cells of the face. $\hat{n} \cdot (\nabla P)^{n+1} \Bigr|_{\substack{f}}$ is computed exactly in the same way as the regular diffusion term by splitting it into direct diffusion and cross diffusion terms (Eq.~(\ref{Eq:Diffusion_2}) with $\phi = P$). The face centered pressure gradient required in the cross diffusion term is estimated by the modified coefficients (Eq.~(\ref{Eq:Face_gradient_modified})). This volume flux estimate satisfying the discrete continuity equation to a specified tolerance is used in the convection term (Eq.~(\ref{Eq:Momentum_frac_step_u*})).
\par The cell centered velocities do not satisfy the discrete continuity equation. They are computed from equation (\ref{Eq:Velocity_correction}) and cell centered pressure gradient is estimated by averaging the face centered gradients.
\subsubsection{Latent Heat Term}
The Gulliver-Scheil equation (\ref{Eq:Solid_fraction}) which relates temperature with solid fraction is a non-linear model. The easiest way to numerically couple this with the energy equation is to model the Gulliver-Scheil equation fully explicitly as a source term. The problem with an explicit approach is that the source term destabilizes the discretized energy equation due to high magnitude of the latent heat coefficient. Thus, we use the source term linearization concept discussed by \citet{patankar1980numerical}. The nonlinear term is split into a linear term and a remainder. The linear term is modeled implicitly and the remainder term is treated explicitly. From numerical stability point of view, the coefficient of the linear term which is added to the diagonal of the matrix should be positive. If the coefficient and the remainder term are functions of the unknown (temperature in this case), the equation has to be solved iteratively till convergence.
\par The latent heat term of the energy equation (\ref{Eq:Energy}) when integrated over time and control volume gives:
\begin{equation}
\int_V \int_t L_f \frac{\partial f_s}{\partial t} dV dt \approx L_f \Delta V \left(f_s ^{m+1}-f_s^{old}\right) 
\label{Eq:Latent_heat_1}
\end{equation}
where, superscripts $old$ and $m$ denote last time-step value and iteration number respectively. The value of solid fraction in the subsequent iteration ($f_s ^{m+1}$) can be estimated from its latest value ($f_s^m$) by a first order Taylor expansion:
\begin{equation}
f_s ^{m+1} \approx f_s^m + \left\{\frac{df_s}{dT}\right\}^m \left[ T_p^{m+1} - T_p^m\right]
\label{Eq:Latent_heat_2}
\end{equation}
Substituting equation (\ref{Eq:Latent_heat_2}) in (\ref{Eq:Latent_heat_1}) gives:
\begin{equation}
\begin{split}
\int_V \int_t L_f \frac{\partial f_s}{\partial t} dV dt &\approx
L_f \Delta V \left(f_s^m + \left\{\frac{df_s}{dT}\right\}^m \left[ T_p^{m+1} - T_p^m\right] - f_s^{old}\right) \\
&= \left[L_f \Delta V \left\{\frac{df_s}{dT}\right\}^m\right] T_p^{m+1} + \left[L_f \Delta V\left(f_s^m - f_s^{old} - \left\{\frac{df_s}{dT}\right\}^m T_p^m\right)\right]\\
&=S_p T_p^{m+1} + S_c
\end{split}
\label{Eq:Latent_heat_3}
\end{equation}
$S_p$ and $S_c$ are functions of last iteration and last time-step values and thus can be computed first. Note that $S_p$ is always negative and when taken to the left hand side of the equation, it becomes positive and is thus added to the diagonal of the linear system matrix. Adding a positive term to the diagonal helps in stabilizing the system and speeds up convergence. Hence, this approach is found to be much better than the fully explicit method.
\begin{figure}[h]
	\centering
	\begin{subfigure}{0.49\textwidth}
		\includegraphics[width=\textwidth]{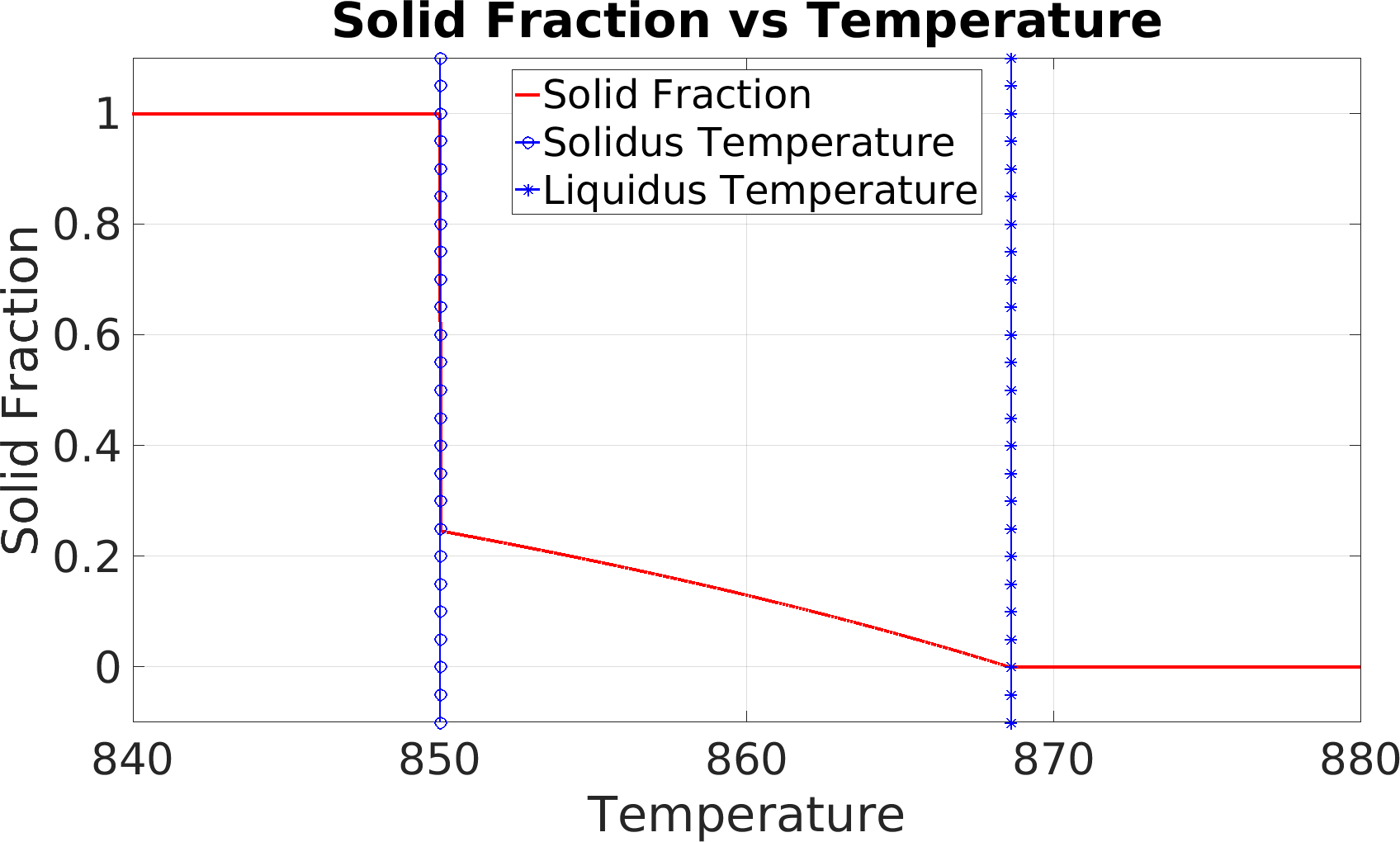}
		\caption{Solid Fraction}
		\label{Fig:solid_fraction}
	\end{subfigure}
	\begin{subfigure}{0.49\textwidth}
		\includegraphics[width=\textwidth]{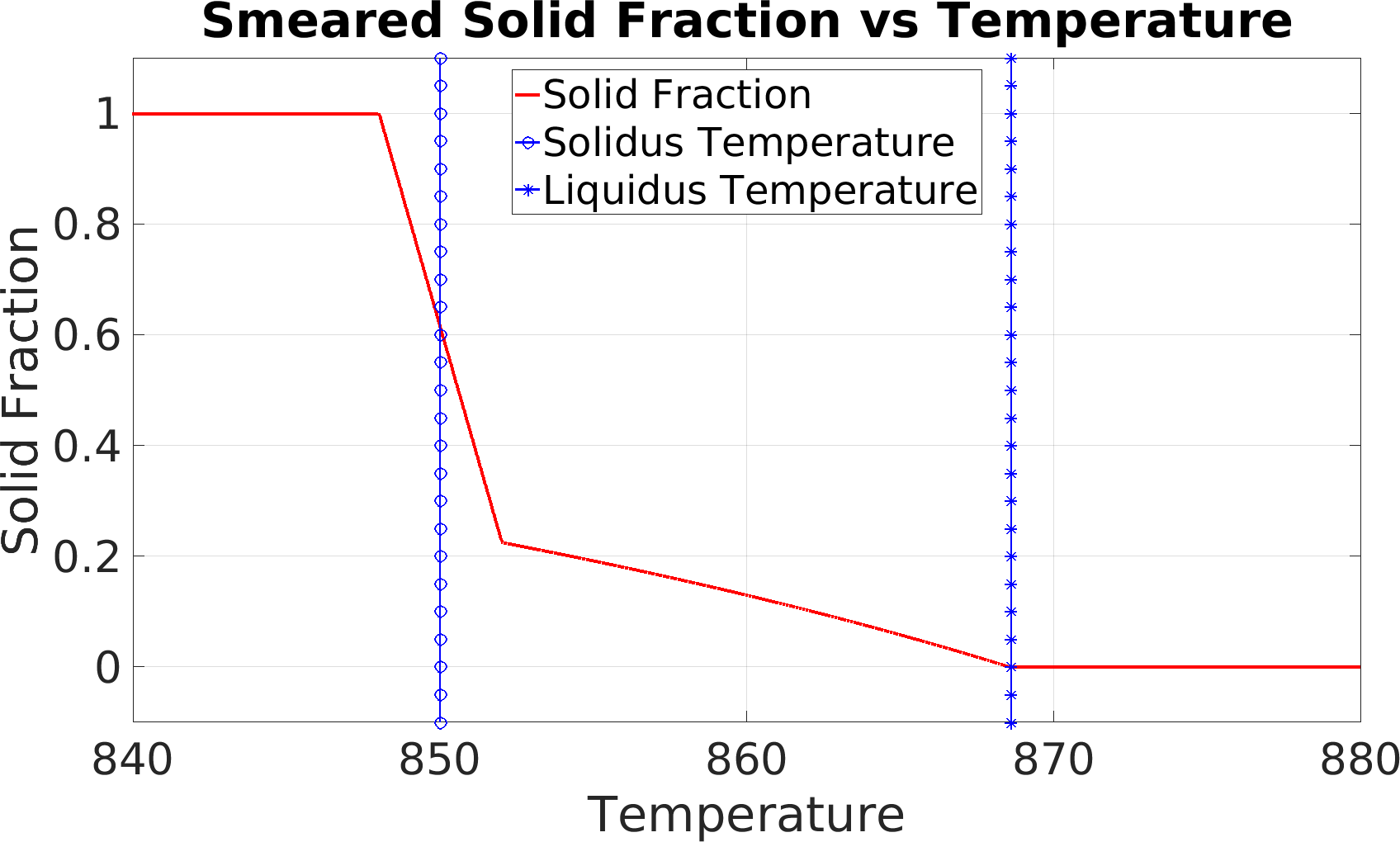}
		\caption{Smeared Solid Fraction}
		\label{Fig:solid_fraction_smeared}
	\end{subfigure}    
	\begin{subfigure}{0.49\textwidth}
		\includegraphics[width=\textwidth]{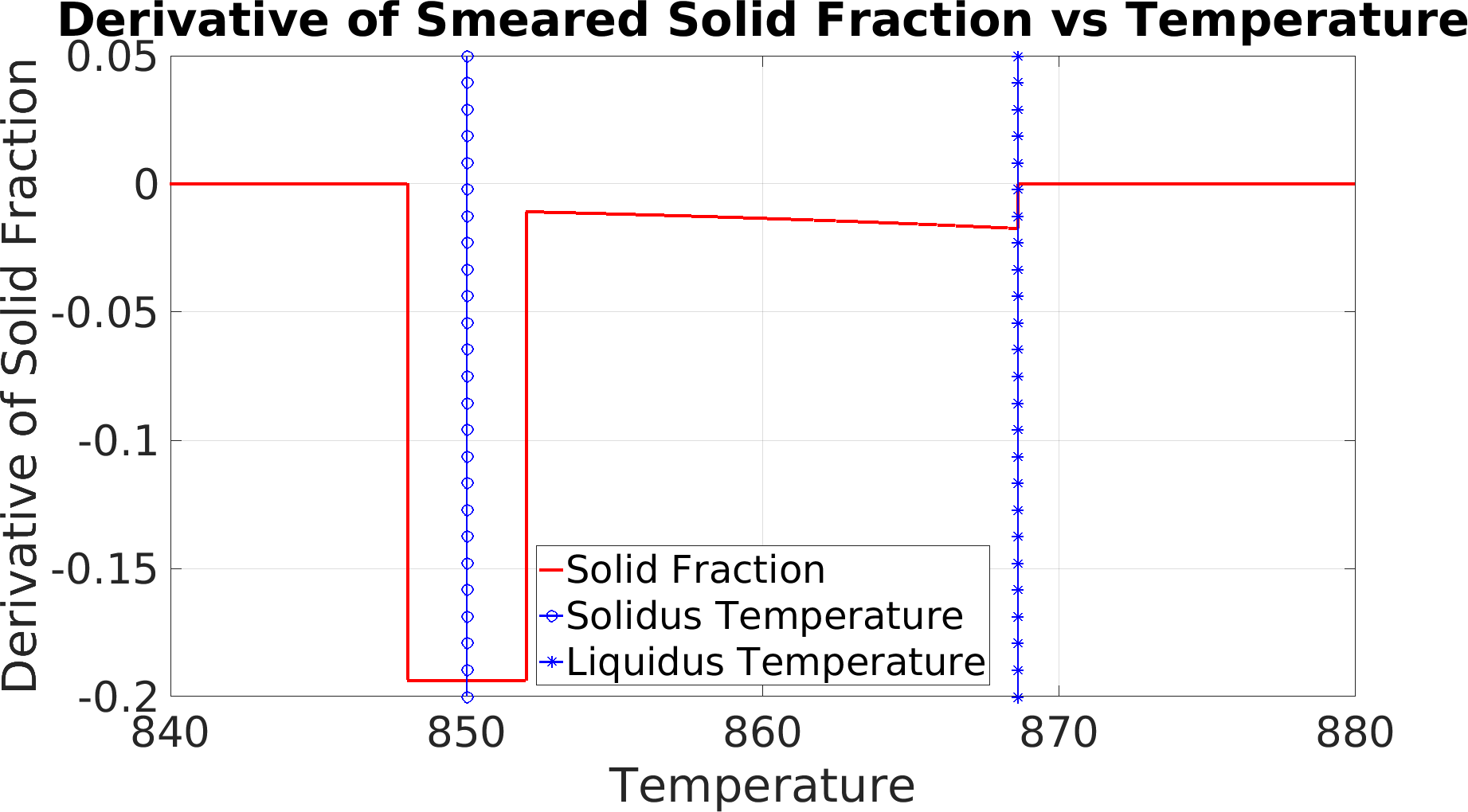}
		\caption{Derivative of Smeared Solid Fraction}	
		\label{Fig:solid_fraction_smeared_derivative}	
	\end{subfigure}
	\caption{}
	\label{Fig:solid_fraction_3_figures}
\end{figure}
\par The Gulliver-Scheil equation (\ref{Eq:Solid_fraction}) plotted in Fig. \ref{Fig:solid_fraction} for a typical aluminum alloy shows that there is a discontinuity at the solidus temperature. Thus, the derivative $\frac{df_s}{dT}$ cannot be computed. To deal with this difficulty, the original equation is modified by smearing the discontinuity near the solidus temperature:
\begin{equation}
f_s(T)=
\begin{cases}
0 & \text{if } T  >T_{liq}\\
1 & \text{if } T  <T_{sol} - T_{\epsilon}\\
\hat{f_s} - (T - T_{sol} - T_{\epsilon}) \left(\frac{1 - \hat{f_s}}{2 T_{\epsilon}}\right)& \text{if } T_{sol} - T_{\epsilon} < T < T_{sol} + T_{\epsilon}\\
1 - \left( \frac{T-T_f}{T_{liq} - T_f}\right)^{\frac{1}{k_p-1}}  & \text{otherwise}
\end{cases}
\label{Eq:Solid_fraction_smeared}
\end{equation}
where, $\hat{f_s} = 1 - \left( \frac{T_{sol} + T_{\epsilon} - T_f}{T_{liq} - T_f}\right)^{\frac{1}{k_p-1}}$ and $T_{\epsilon}$ is the width of linear smear which can be set to a reasonable value like 2 K. Thus, the derivative can be computed analytically. Figures~\ref{Fig:solid_fraction_smeared} and \ref{Fig:solid_fraction_smeared_derivative} plot the modified solid fraction relation and its derivative respectively.
\par The overall iterative procedure to obtain the variables at the new time step from values at the old time step can be summarized as:
\begin{enumerate}
	\item Initialize: $T_p^0 = T_p^{old}$ and $f_s^0 = f_s^{old}$
	\item Compute $S_p$ and $S_c$ using last iteration values ($T_p^m$ and $f_s^m$) by equation (\ref{Eq:Latent_heat_3}) and solve the linear system of equations to estimate next iteration value $T_p^{m+1}$
	\item Update the solid fraction: $f_s ^{m+1} = (1-\lambda) f_s^m + \lambda f_s(T_p^{m+1})$ where, $0<\lambda \leq 1$ is an under relaxation parameter
	\item Estimate the relative change between the successive iteration values of temperature and solid fraction
\end{enumerate}
Repeat steps $2-4$ until the relative change drops below a desired threshold. For the aluminum alloy used here, it is found that under relaxation is not required i.e., $\lambda = 1$ and the solution converges in $5-10$ iterations.
\subsection{Solver for Linear Systems}
\par Typical die cast geometries have high aspect ratios i.e., thin cross sections compared to the lateral dimensions. It is found that single grid iterative solvers for the elliptic pressure Poisson equation converge slowly for such geometries. Hence, in this work a multigrid solver is used. The central idea of a multigrid solver is to solve the equations on multiple coarse grids and couple the corrections from all the grids through prolongation and relaxation. The high frequency component of the residual converges fast on the fine grids while the coarse grids are used to accelerate convergence of the low frequency residual. Thus, the coarse grid solutions are used to accelerate the convergence while maintaining the discretization accuracy of the solution at the finest level.
\par Geometric multigrid is a technique in which multiple levels of grids are generated physically and the matrix vector system is estimated by discretizing the governing equations at each level. The main benefit of this approach is that the matrices at all the levels are obtained directly from the governing equations and thus, good convergence is observed. The main drawback is that generating coarse grids for a complex geometry with unstructured elements is non-trivial. Algebraic multigrid (AMG) tries to address this problem by coarsening the matrix using heuristics based algorithms. This is a black box approach which does not need any physical grids at coarse levels.
\par The BoomerAMG routine along with Krylov solvers of the open source library HYPRE \cite{HYPRE_website} developed at the Lawrence Livermore National Laboratory is used in our work. The AMG solver is found to solve the modified momentum equations (\ref{Eq:Momentum_frac_step_u*}) and the energy equation (\ref{Eq:Energy}) without difficulty. However, some consistency issues have been observed during the solution of the pressure Poisson equation (\ref{Eq:Pressure_poisson}) when extensive solidification happens. In complex geometries, as solidification proceeds, there can be disjoint pockets of metal which are yet to solidify. As these pockets are far enough from each other, they are decoupled numerically in the discrete reduced pressure Poisson equation (section \ref{Sec:Pressure Poisson Equation}). Hence, in some cases, there are zeros on the diagonals of the coarsest grid level generated by the heuristic coarsening algorithm of BoomerAMG, which creates problems in the convergence of the solution. Thus, the AMG solver has been coupled with a single grid BiCGSTAB solver, also from HYPRE, to be used when AMG is unable to solve the pressure Poisson equation. OpenCast has conditions programmed which automatically switch to BiCGSTAB when the AMG solver diverges. Note that this issue arises towards the end of the simulation when only a few cells are liquid (for instance 20\%). Thus, the reduced system is much smaller in size compared to the original problem and a single grid solver is reasonably well convergent.
\subsection{Grain Growth and Mechanical Properties Models}
Grain size and Secondary Dendrite Arm Spacing (SDAS) are two important parameters used to characterize the microstructure. OpenCast uses empirical relations from the literature for estimation of microstructure parameters and mechanical properties such as yield strength.
\par SDAS is predicted based on the empirical relationship proposed by \citet{backer2007microporosity}, 
\begin{equation} 
SDAS = \lambda_2 =A_\lambda\left(\frac{\partial T}{\partial t}\right)^{B_\lambda} \hspace{10pt} [in~\mu m].
\label{Eq:SDAS}
\end{equation}
The model parameters $A_\lambda$ and $B_\lambda$ are chosen to be 39.4 and -0.317, respectively based on the model for microstructure in aluminum alloys \cite{backer2007microporosity}. Material behavior of the die cast alloy is predicted in terms of 2\% yield strength ($\sigma_{0.2}$) using empirical relationship proposed by \citet{okayasu2015precise}.
\begin{equation} 
\sigma_{0.2} = A_\sigma \lambda_2 ^{-1/2} + B_\sigma
\label{Eq:YS}
\end{equation}
Here, $\sigma_{0.2}$ is in MPa, $\lambda_2$ (SDAS) is in $\mu m$, $A_\sigma=59.0$ and $B_\sigma=120.3$~\cite{okayasu2015precise}.
\par Grain size estimation is based on the work of \citet{greer2000modelling}. The grain growth rate is given by:
\begin{equation} 
\frac{dr}{dt} = \frac{\lambda_s^2 D_s}{2 r}
\label{Eq:grain_growth_rate}
\end{equation}
where, $r$ is the grain size, $D_s$ is the solute diffusion coefficient in the liquid and $t$ is the time. The parameter $\lambda_s$ is obtained by invariant size approximation:
\begin{equation} 
\lambda_s = \frac{-S}{2 \pi^{0.5}} + \left(\frac{S^2}{4 \pi} - S\right)^{0.5}
\label{Eq:grain_growth_lambda}
\end{equation}
$S$ is given by 
\begin{equation} 
S = \frac{2(C_s - C_0)}{C_s - C_l}
\label{Eq:grain_growth_S}
\end{equation}
where, $C_l = C_0 (1-f_s)^{(k_p-1)}$ is solute content in the liquid, $C_s = k_p C_l$ is solute content in the solid at the solid-liquid interface and $C_0$ is the nominal solute concentration. Thus, using the prescribed partition coefficient ($k_p$) and estimated solid fraction ($f_s$), equations (\ref{Eq:grain_growth_rate})--(\ref{Eq:grain_growth_S}) are solved to get the final grain size.
\subsection{Parameter Uncertainty Quantification}\label{Sec:Parameter Uncertainty Quantification}
The final product quality in die casting is influenced by process parameters like alloy material properties, interface conditions at the mold, thermal boundary conditions etc. Due to the complexity of the process, accurate measurement and control of these parameters is difficult. Conventional deterministic simulations alone are not sufficient to estimate the effect of any stochastic variations on the product quality, thus parameter uncertainty quantification is important. From modeling point of view, parameter uncertainty quantification is a set of stochastic partial differential equations with variation in initial conditions, coefficients and boundary conditions. Stochastic variables are considered as dimensions of the problem in addition to space and time.
\par To estimate the relation between stochastic process parameters and output parameters, various methods have been proposed in the literature. A popular method is to use a linear combination of polynomial basis functions in the stochastic dimension to expand the output variables. Since orthogonality helps in convergence, orthogonal polynomials are used as basis functions. Wiener's polynomial chaos \cite{wiener1938homogeneous} with Askey family of orthogonal polynomials leads to optimal convergence of the error. \citet{xiu2002wiener} determined which polynomial leads to exponential convergence depending on the underlying probability distribution function that the stochastic variable follows. For instance, the Hermite polynomials are orthogonal with respect to the standard normal distribution as the weighting function. Hence, it is recommended to use Hermite polynomials as basis functions if the stochastic variable follows normal distribution. Using generalized polynomial chaos expansion developed by \citet{xiu2002wiener}, a second order random field ($w$) can be expanded by polynomial basis (Eq.~(\ref{Eq:Polynomial_chaos})). For all practical purposes, the series is truncated to order $n$.
\begin{equation} 
w(\bm{x},\bm{\xi}(\theta))= \sum_{i=0}^{\infty} w_i(\bm{x})\bm{\Psi_i}(\bm{\xi}(\theta))\approx\sum_{i=0}^{n} w_i(\bm{x})\bm{\Psi_i}(\bm{\xi}(\theta))
\label{Eq:Polynomial_chaos}
\end{equation}
\par In this work, stochastic collocation method is used to estimate the deterministic coefficients of the polynomial chaos expansion. Collocation is a non-intrusive method and thus, modification of deterministic software is not necessary. It acts as a wrapper around existing deterministic software. The deterministic simulation is run at $M$ sample points ($\bm{\xi}^m$) and a condition $w(\bm{x},\bm{\xi}^m)=w_{sim}(\bm{x},\bm{\xi}^m)$ is imposed. The right hand side comes from each deterministic simulation and left hand side from polynomial expansion. This gives $M$ constraints written in a matrix vector form \cite{smith2013uncertainty}. $M>n+1$ ensures that the Vandermonde system (Eq.~(\ref{Eq:Collocation_vandermonde})) is overdetermined. Solving in the least-squares sense gives $\begin{bmatrix} w_0\left(\bm{x}\right)&  \cdots & w_n\left(\bm{x}\right) \end{bmatrix} ^T$.
\begin{equation} 
\renewcommand\arraystretch{1.8}
\begin{bmatrix} \bm{\Psi_0}\left(\bm{\xi}^1\right) & \cdots & \bm{\Psi_n}\left(\bm{\xi}^1\right)  \\ \vdots  && \vdots \\ \bm{\Psi_0}\left(\bm{\xi}^M\right) & \cdots & \bm{\Psi_n}\left(\bm{\xi}^M\right) \end{bmatrix}
\begin{bmatrix}
w_0(\bm{x})\\
\vdots \\
w_n(\bm{x})
\end{bmatrix} = 
\begin{bmatrix}
w_{sim}\left(\bm{x},\bm{\xi}^1\right)\\
\vdots \\
w_{sim}\left(\bm{x},\bm{\xi}^M\right)
\end{bmatrix} 
\label{Eq:Collocation_vandermonde}
\end{equation}
\par Sample points ($\bm{\xi}^m$) have to be chosen wisely for successive implementation of stochastic collocation method. For instance, uniformly distributed samples can lead to highly oscillatory basis functions and thus poor convergence. Thus, for one dimensional stochastic problems \cite{smith2013uncertainty}, it is popular to choose roots of the basis orthogonal polynomial as sample points. For multiple stochastic dimensions, a simple idea is to use tensor product of the single dimensional samples. The problem with tensor products is that the number of samples grows exponentially with stochastic dimensions. Each sample corresponds to a deterministic simulation and thus the computational expense grows exponentially. \citet{smolyak1963quadrature} came up with an algorithm to reduce number of samples in multi-dimensional space maintaining the accuracy of the interpolation. In this research, sparse grid nodes are taken from the work of \citet{heiss2008likelihood}. A MATLAB based tool UQLab developed by \citet{marelli2014uqlab} is used for post processing the simulation outputs to estimate polynomial chaos coefficients and generate the response surfaces.
\par The strategy described above is quite general and can be applied to any numerical solution framework. In this case, the output variables ($w$ of Eq.~\ref{Eq:Polynomial_chaos}) are temperature and microstructure parameters like grain size and dendritic arm spacing. The stochastic input parameters are boundary conditions, initial conditions and alloy and material properties. The polynomial chaos expansion with stochastic collocation is combined with the deterministic computational fluid flow and heat transfer solver to estimate the sensitivity and uncertainty propagation.
\section{Validation with Uncertainty}
Due to the complexity of the die casting process, controlled experiments with accurate temperature measurements inside the casting during solidification are difficult. To validate our code, we have therefore used the experimental results of ingot solidification made of Sn-Bi alloy reported by \citet{quillet2007benchmark}. In their experiments, a $50 \times 60 \times 10$ mm ingot is solidified by cooling one of the $60 \times 10$ mm faces at a constant rate (5 K/min) and the remaining 5 faces are thermally insulated. Twenty five thermocouples on the ingot were used to measure temperatures during solidification. Due to the low cooling rate, natural convection velocities were significant and their effect can be seen from the curved contours of temperature. Thus, this validation study covers all the aspects of the software including heat transfer, natural convection and solidification.
\par
\begin{figure}[h]
	\centering
	\begin{subfigure}{0.51\textwidth}
		\includegraphics[width=\textwidth]{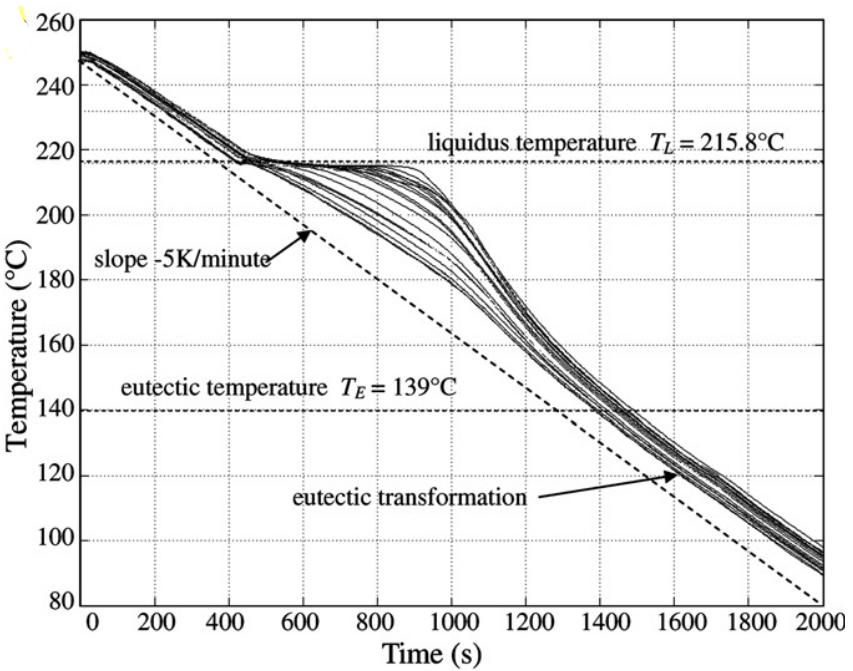}
		\caption{Experiment \cite{quillet2007benchmark}}
		\label{Fig:exp_Temp_time_cooling_rate_5}
	\end{subfigure}	
	\begin{subfigure}{0.49\textwidth}
		\includegraphics[width=\textwidth]{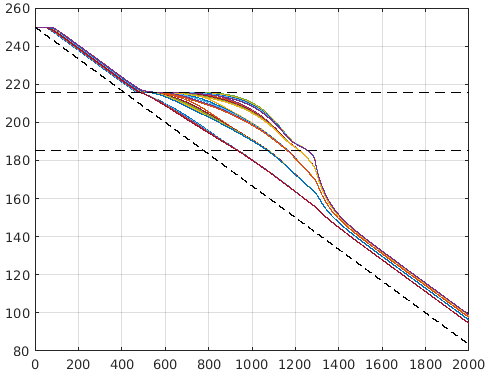}
		\caption{Simulation with Constant Thermal Conductivity}	
		\label{Fig:sim_Temp_time_cooling_rate_5_cons_k_30_Teps_5}	
	\end{subfigure}
	\begin{subfigure}{0.49\textwidth}
		\includegraphics[width=\textwidth]{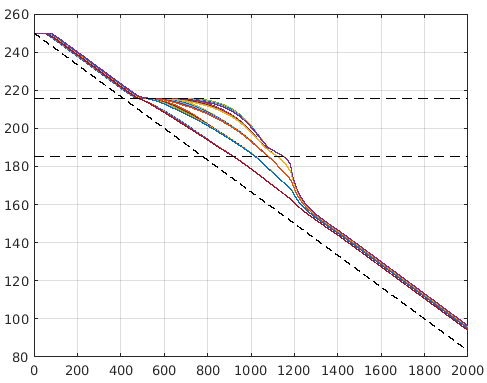}
		\caption{Simulation with Variable Thermal Conductivity}
		\label{Fig:sim_Temp_time_cooling_rate_5_var_k_Teps_5}
	\end{subfigure}    
	\caption{Temperature-Time Plot for Cooling Rate 5 K/min}
	\label{Fig:Temp_time_cooling_rate_5}
\end{figure}

 As a first step, validation was attempted without incorporating the effects of uncertainty. Figure~\ref{Fig:exp_Temp_time_cooling_rate_5} plots the experimental measurements of temperature from the 25 thermocouples as a function of time beginning from solidification. Figure~\ref{Fig:sim_Temp_time_cooling_rate_5_cons_k_30_Teps_5} plots the temperatures predicted by OpenCast simulations with a constant thermal conductivity. Subsequently, we included the temperature varying thermal conductivity which changes the temperature time curves to Fig.~\ref{Fig:sim_Temp_time_cooling_rate_5_var_k_Teps_5}. This observation shows that the temperature dependency of thermal conductivity has an appreciable effect on the temperature evolution.
\par Local values of the thermo-physical properties such as thermal conductivity and density depend on grain structure. At the length scale of current simulation, it is not possible to predict the grain structure and hence the properties from first principles. Thus, there can be some uncertainty in the input properties. Also, when the alloy solidifies and cools, there is a thermal contraction. This creates a gap between the mold wall and the casting and adds thermal contact resistance, reducing the amount of heat extracted. Further, the experimental measurements are also subject to sensor noise. For example, the thermocouples used by \citet{quillet2007benchmark} have a reported accuracy of 1K, thereby adding another uncertainty in the predictions. Thus, the validation must account for these uncertainties. We have therefore included the stochastic variation in the boundary conditions and material properties and propagated them using the stochastic collocation method. A confidence interval is estimated about the mean output parameters and the experimental and numerical predictions are compared within the confidence interval.
\par Uncertainty is therefore added to wall temperature, latent heat, density and thermal conductivity. Since it is difficult to estimate the thermal contact resistance, the wall temperature is specified as a Dirichlet boundary condition by adding an offset to the cooling rate of 5 K/min. The offset is estimated from the experimental temperature plot (Fig.~\ref{Fig:exp_Temp_time_cooling_rate_5}). In order to take into account the errors in the temperature measurement, an uncertainty is added on top of the offset. All the four input stochastic parameters are assumed to follow a normal distribution with mean ($\mu$) and standard deviation ($\sigma$) as follows:
\begin{enumerate}
	\item Wall temperature offset: $\mu = 0^0$ C, $\sigma = 0.5^0$ C
	\item Latent heat: $\mu = 60000$ J/kg, $\sigma = 2\% \mu$ 
	\item Density: $\mu = 7300$ kg/m$^3$, $\sigma = 2\% \mu$ 
	\item Temperature dependent thermal conductivity: $\mu = [61.282,57.42,30.1,37.7]$ W/mK, $\sigma = [2.5,2.5,2.5,2.5]$ W/mK at temperature $[273.2, 373.2, 573.2, 973.2]$ K
\end{enumerate}
All the thermo-physical properties are estimated as a weighted average of individual properties of Sn and Bi taken from the online version of \citet{kaye1921tables}. These uncertainties are propagated using the stochastic collocation method described in section \ref{Sec:Parameter Uncertainty Quantification}.
\par In order to make a comparison with the experimental data, the experimental temperature-time data plot (Fig.~\ref{Fig:exp_Temp_time_cooling_rate_5}) is digitized. Since it is difficult to distinguish between the 25 temperatures, the thermocouple with the highest temperature is used for validation. Stochastic collocation is done with three accuracy levels of sample points in order to study its convergence. For estimating the interpolation error, 60 Latin Hypercube samples are used as test points. Deterministic simulations and polynomial chaos give two independent estimates of the same output parameter at the test points. The non-dimensional error is defined as the root mean square of difference between these two estimates divided by the maximum value of the parameter. First column of Table~\ref{Table:Stochastic Collocation Convergence Analysis} denotes the accuracy level of sample points. Accuracy level $l$ integrates polynomials upto order $2l-1$ exactly \cite{heiss2008likelihood}. Second column is the number of sample points in two dimensional Smolyak sparse grid i.e., the number of deterministic simulations required ($M$ in Eq.~(\ref{Eq:Collocation_vandermonde})). The last column lists the non-dimensional RMS error in computation of the temperature. It can be seen that the error is of order $10^{-4}$ for all the accuracy levels and it decreases with increasing level thus, showing convergence. Hence, level 6 is used for validation.
\begin{table}[h]
	\begin{center}
		\begin{tabular}{|c|c|c|c|}
			\hline
			Accuracy Level&\# Sample Pts. &Max Temp. RMS Error\\
			\hline 
			4&137 & 7.27E-4\\
			5&385 & 4.11E-4\\
			6&953 & 2.60E-4\\
			\hline
		\end{tabular}
	\caption{Stochastic Collocation Convergence Analysis}
	\label{Table:Stochastic Collocation Convergence Analysis}
	\end{center}
\end{table}
\par The Polynomial-Chaos-Kriging (PCK) module of UQLab \cite{marelli2014uqlab} is used to estimate the output (temperature--time history) as a function of the four input stochastic parameters. Kriging is a Gaussian process modeling tool based on stochastic interpolation algorithm to relate the inputs to outputs. PCK combines the polynomial chaos and kriging methods in a way that is more efficient than the individual methods.
\begin{figure}[h]
	\centering
	\includegraphics[width=0.75\textwidth]{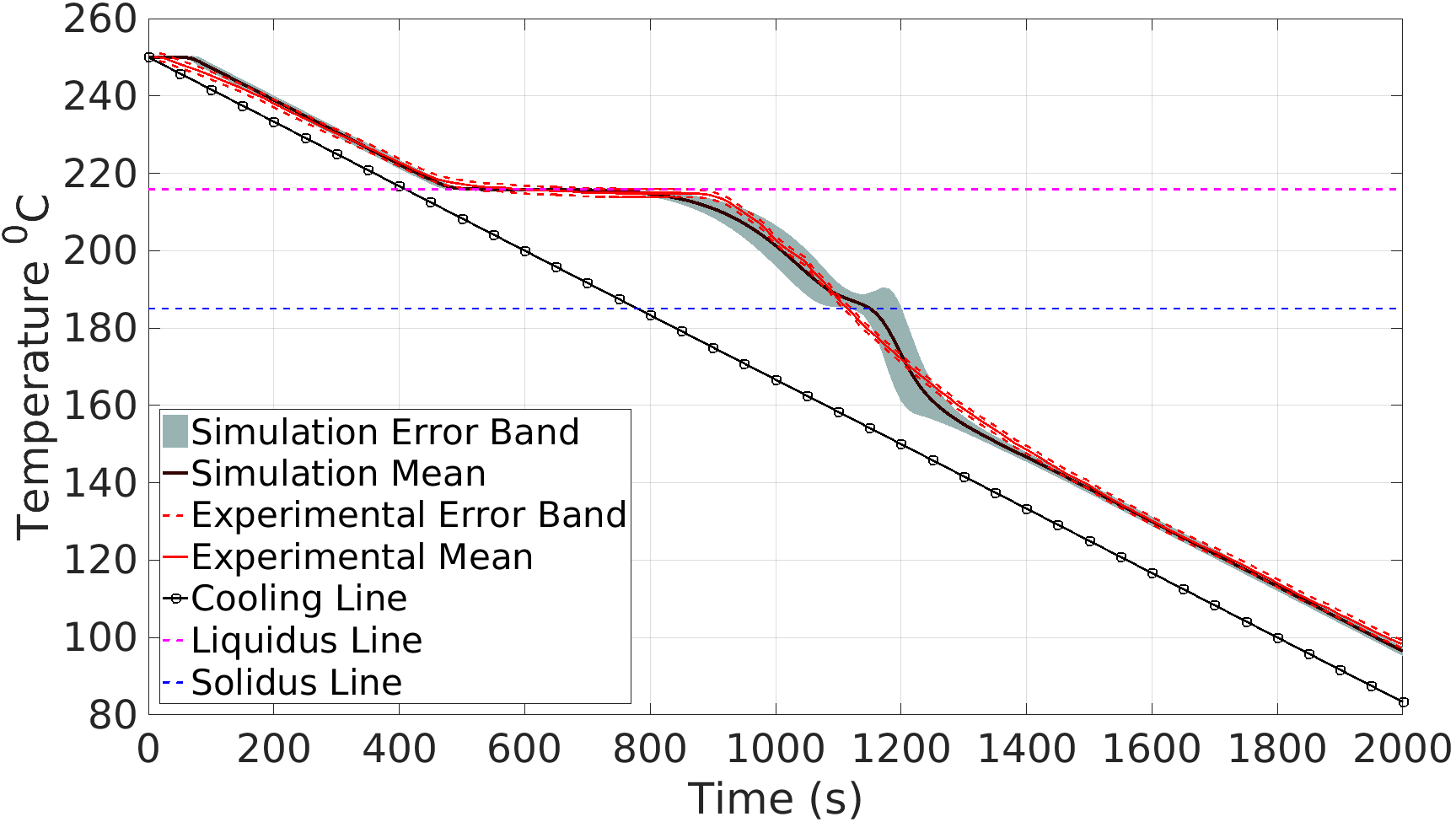}
	\caption{Experimental and Simulated Temperatures with Error Bounds}
	\label{Fig:Max_Temp_Exp_Sim_Error_Bands}
\end{figure}
\par Figure~\ref{Fig:Max_Temp_Exp_Sim_Error_Bands} plots the maximum temperatures (with error bands) from OpenCast simulations and experiments \cite{quillet2007benchmark} together with cooling, solidus and liquidus lines. Before the solidification begins, the temperature drops at a rate similar to the cooling line (5 K/min). Thus, the temperature curve and cooling line are parallel. At the liquidus line, a plateau is observed in the temperature curve. This is due to the latent heat release inside the plate. Near the solidus temperature, it can be seen that the numerical simulation shows a small kink in the temperature curve whereas, the experimental curve drops smoothly. This can be attributed to the numerical smearing of the solid fraction$-$temperature relation (Eq.~(\ref{Eq:Solid_fraction_smeared})). There is a sudden release of latent heat due to sharp change in solid fraction near the solidus temperature (Fig.~\ref{Fig:solid_fraction_3_figures}). It is found that the kink reduces if the width of linear smear ($T_{\epsilon}$ in Eq.~(\ref{Eq:Solid_fraction_smeared})) is increased. The $2 \sigma$ confidence interval about the mean is plotted as the simulation error band. Since the accuracy of experimental temperature measurement is 1 K \cite{quillet2007benchmark}, an experimental error band is plotted as $\pm 1$ K about the mean value. It can be seen that the experimental and numerical estimates are in good agreement thus validating the computational tool. This validation study also emphasizes the utility of uncertainty quantification as previously there was a mismatch between experiments and computations when no uncertainties were included.
\section{Deterministic Solidification Results of Realistic Geometries}
\begin{figure}[h]
	\begin{subfigure}{0.54\textwidth}
		\includegraphics[width=\textwidth]{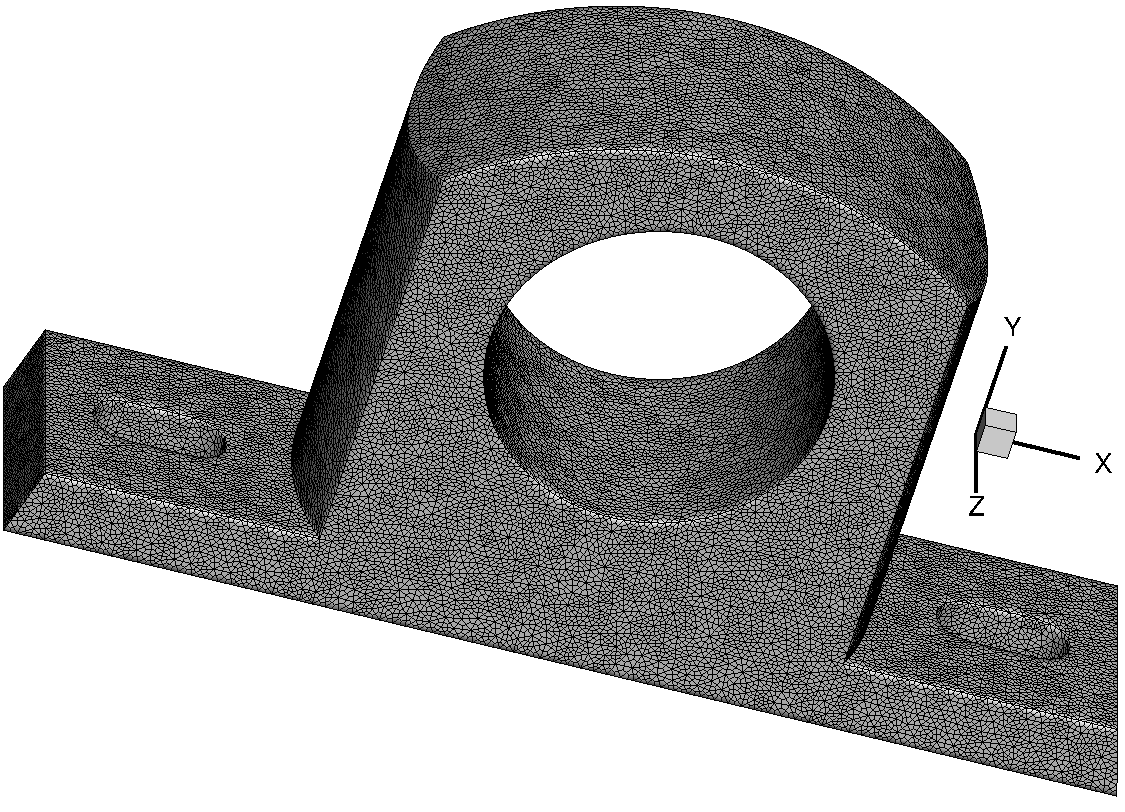}
		\caption{Clamp: 16.5 cm x 9 cm x 3.7 cm}
		\label{Fig:clamp mesh}
	\end{subfigure}	
	\begin{subfigure}{0.45\textwidth}
		\includegraphics[width=\textwidth]{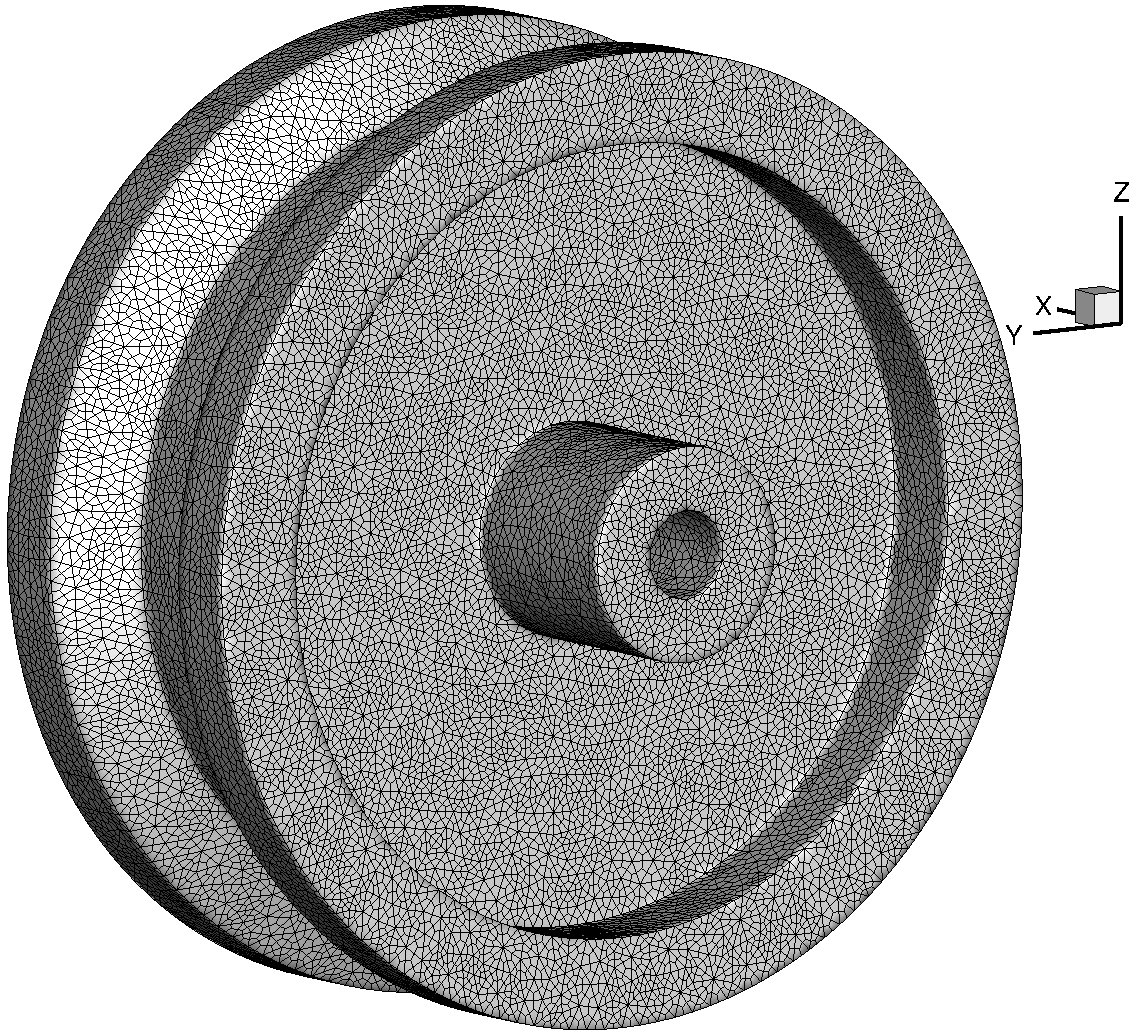}
		\caption{Pulley: 13 cm x 27 cm x 27 cm}	
		\label{Fig:pulley mesh}	
	\end{subfigure}
	\caption{Geometries and the Hexahedral Meshs}
	\label{Fig:mesh}
\end{figure}
In order to demonstrate the utility of the developed software, some complex die casting geometries were simulated. Many die casting geometries have thin cross sections and solidify quickly in a few seconds. Thus, natural convection velocities are negligible and may have little effect on solidification. Thus, simulating without natural convection may be acceptable as it saves significant computational effort. Hence, in all the simulations of this section, natural convection was neglected.
\begin{figure}[H]
	\begin{subfigure}{0.23\textwidth}
		\includegraphics[width=\textwidth]{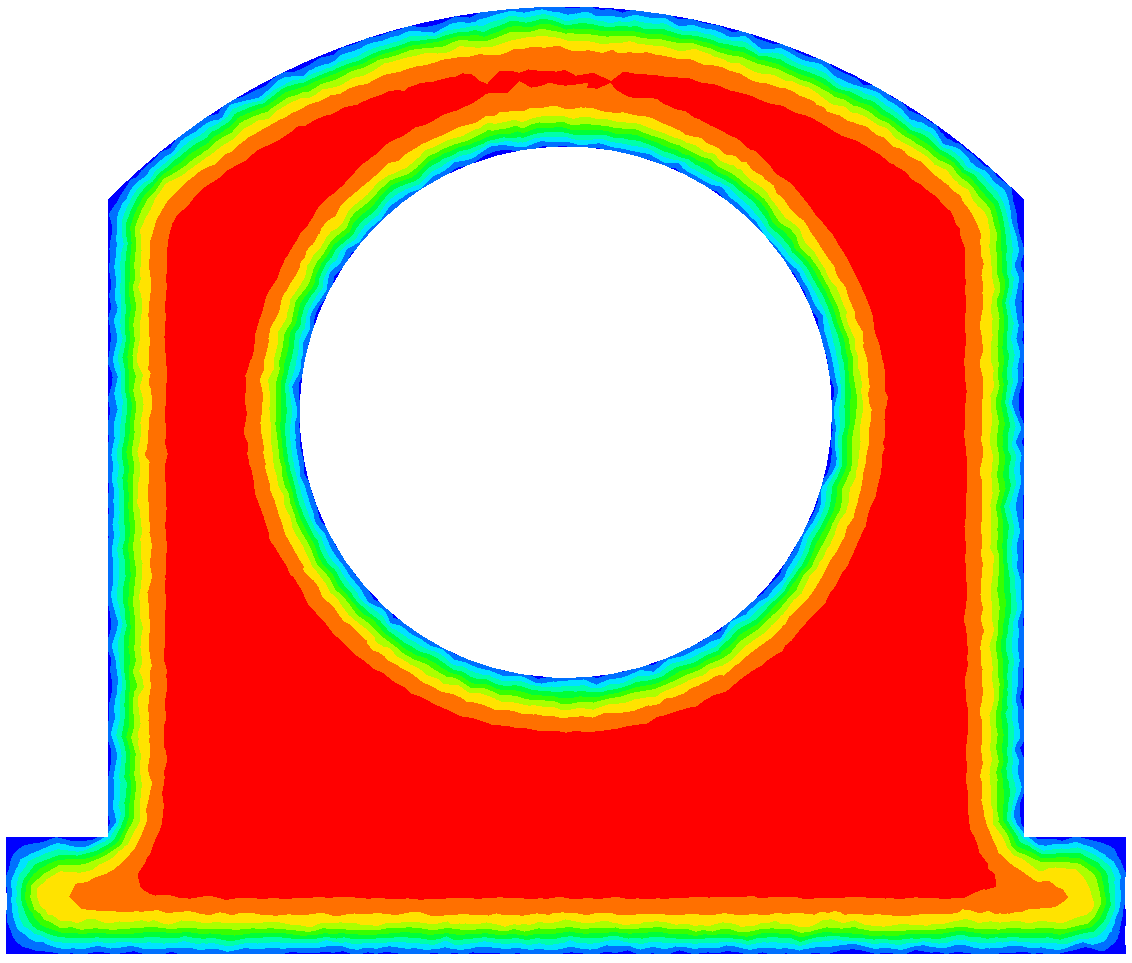}
		\caption{Time: 0.119 s}
	\end{subfigure}	
	\begin{subfigure}{0.23\textwidth}
		\includegraphics[width=\textwidth]{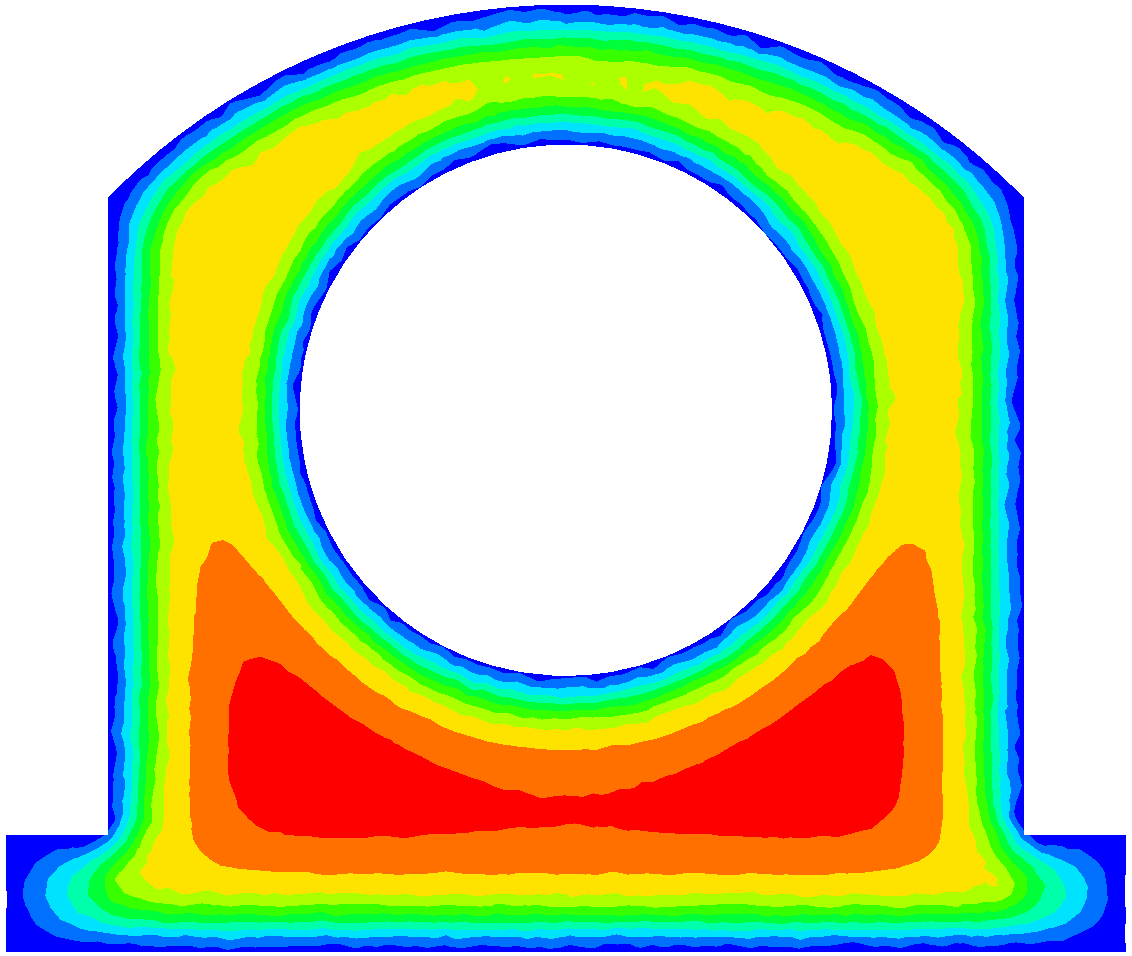}
		\caption{Time: 0.480 s}
	\end{subfigure}
	\begin{subfigure}{0.23\textwidth}
		\includegraphics[width=\textwidth]{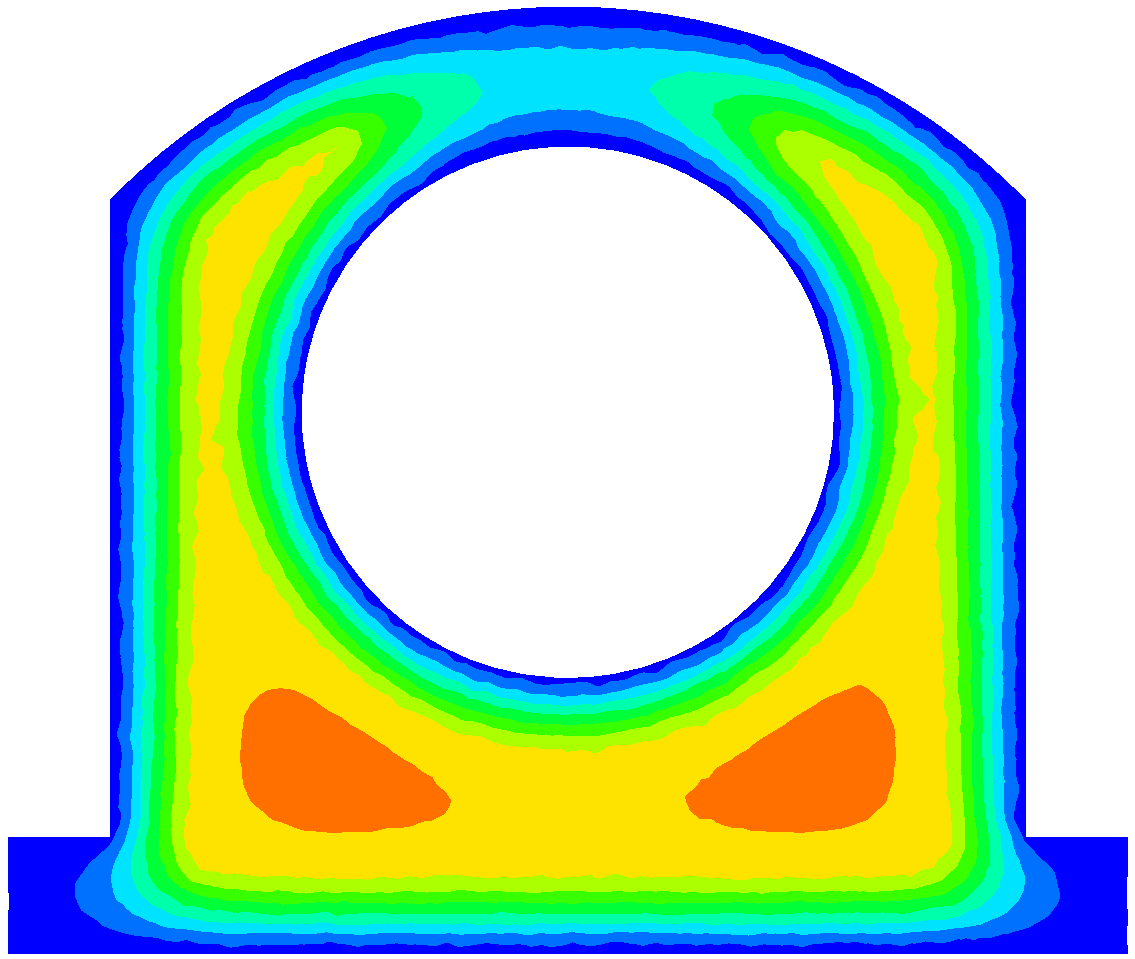}
		\caption{Time: 0.848 s}
	\end{subfigure}
	\begin{subfigure}{0.23\textwidth}
		\includegraphics[width=\textwidth]{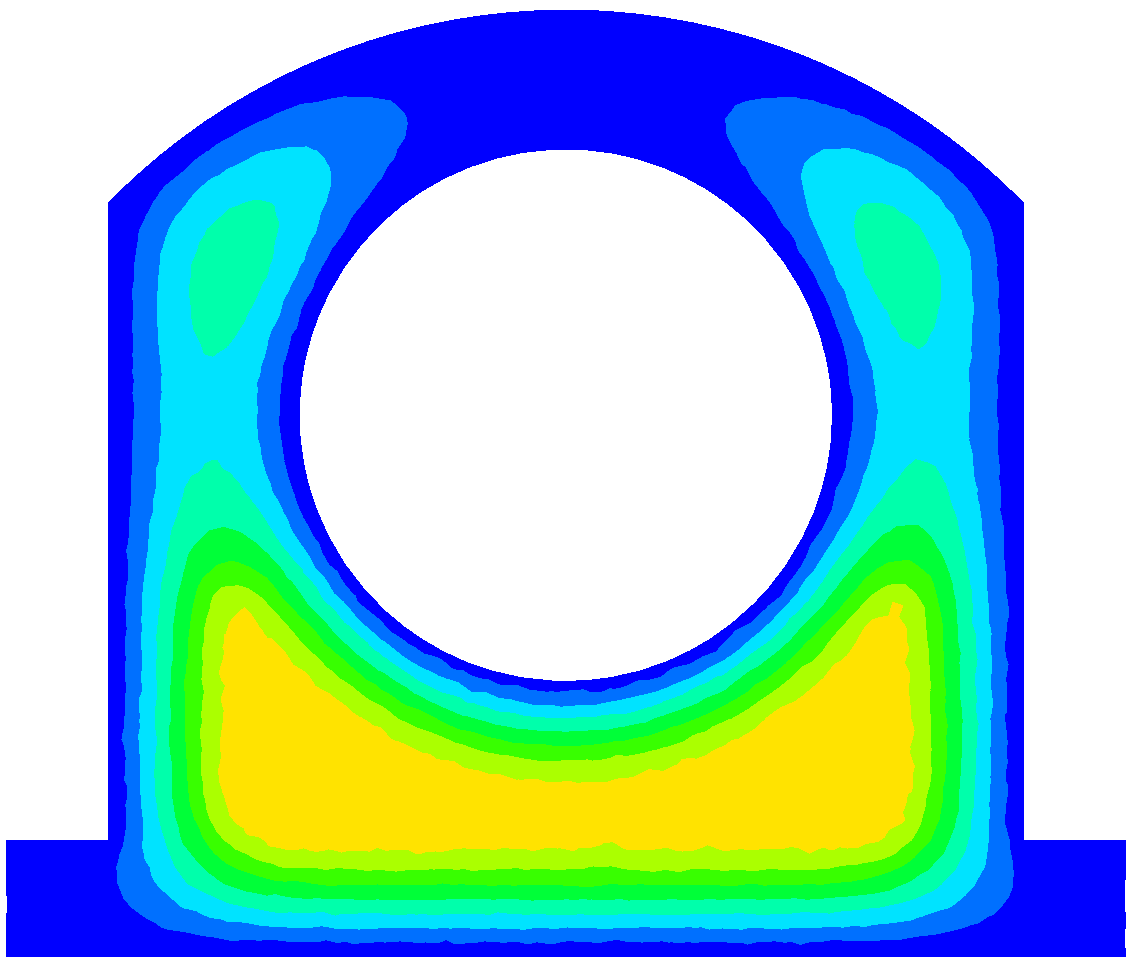}
		\caption{Time: 1.42 s}
	\end{subfigure}	
	\begin{subfigure}{0.05\textwidth}
		\includegraphics[width=\textwidth]{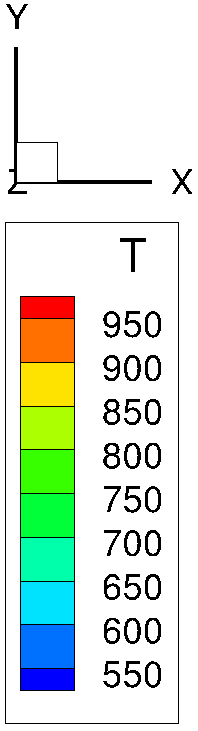}
	\end{subfigure}	
	\caption{Clamp Temperature Contours (K) (Solidification Time: 2.40 s)}
	\label{Fig:clamp_T}
\end{figure}
\par Figure~\ref{Fig:mesh} shows the hexahedral grids of two selected geometries. Both the clamp and the pulley have approximately 300,000 control volumes each. Initially, the mould is filled with a liquid alloy at 1000 K and the walls are held at 500 K. Aluminum alloys with around 8-12\% Silicon content are popular in die casting. Thus, for this simulation, an Al-10\%Si binary alloy is chosen. The liquidus, solidus and freezing temperatures are taken from the phase diagram.
\begin{figure}[H]
	\begin{subfigure}{0.23\textwidth}
		\includegraphics[width=\textwidth]{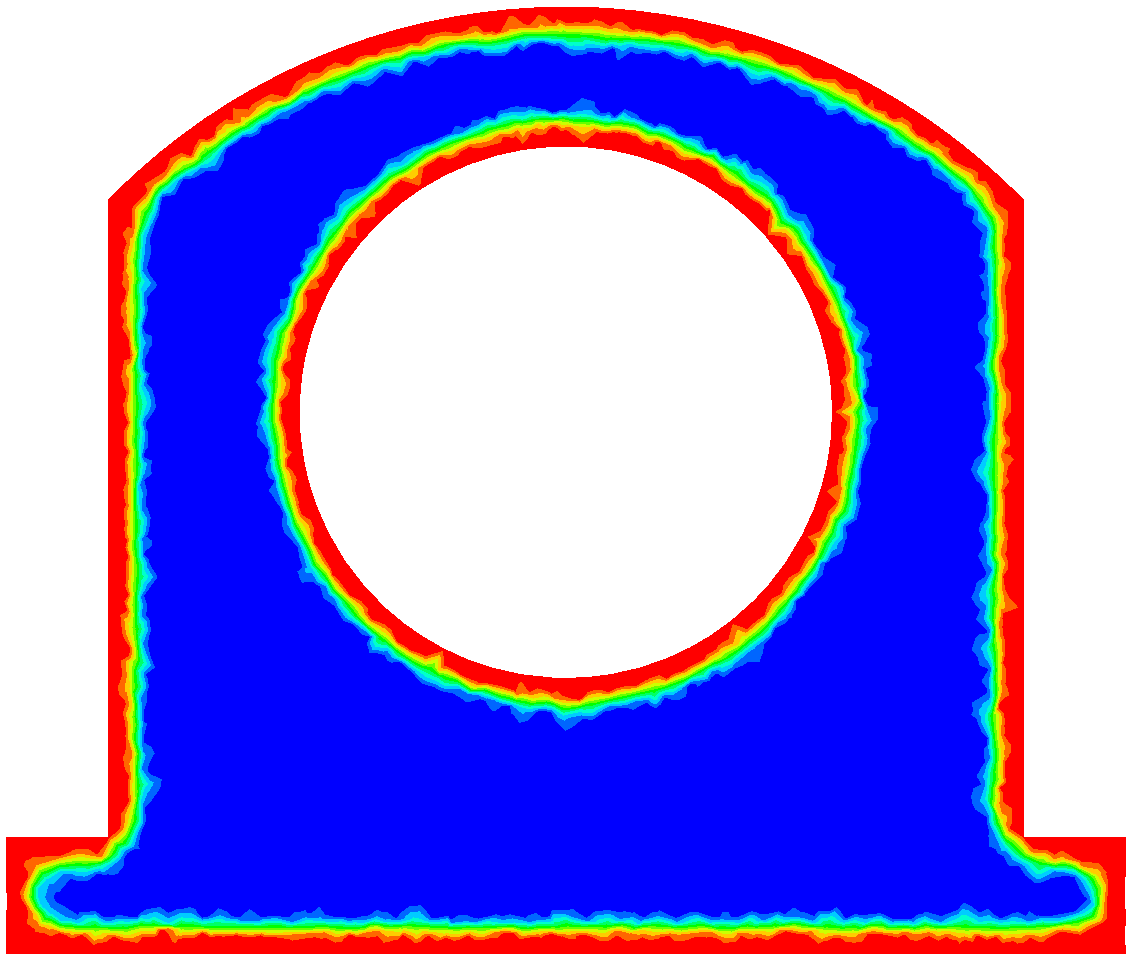}
		\caption{Time: 0.119 s}
	\end{subfigure}	
	\begin{subfigure}{0.23\textwidth}
		\includegraphics[width=\textwidth]{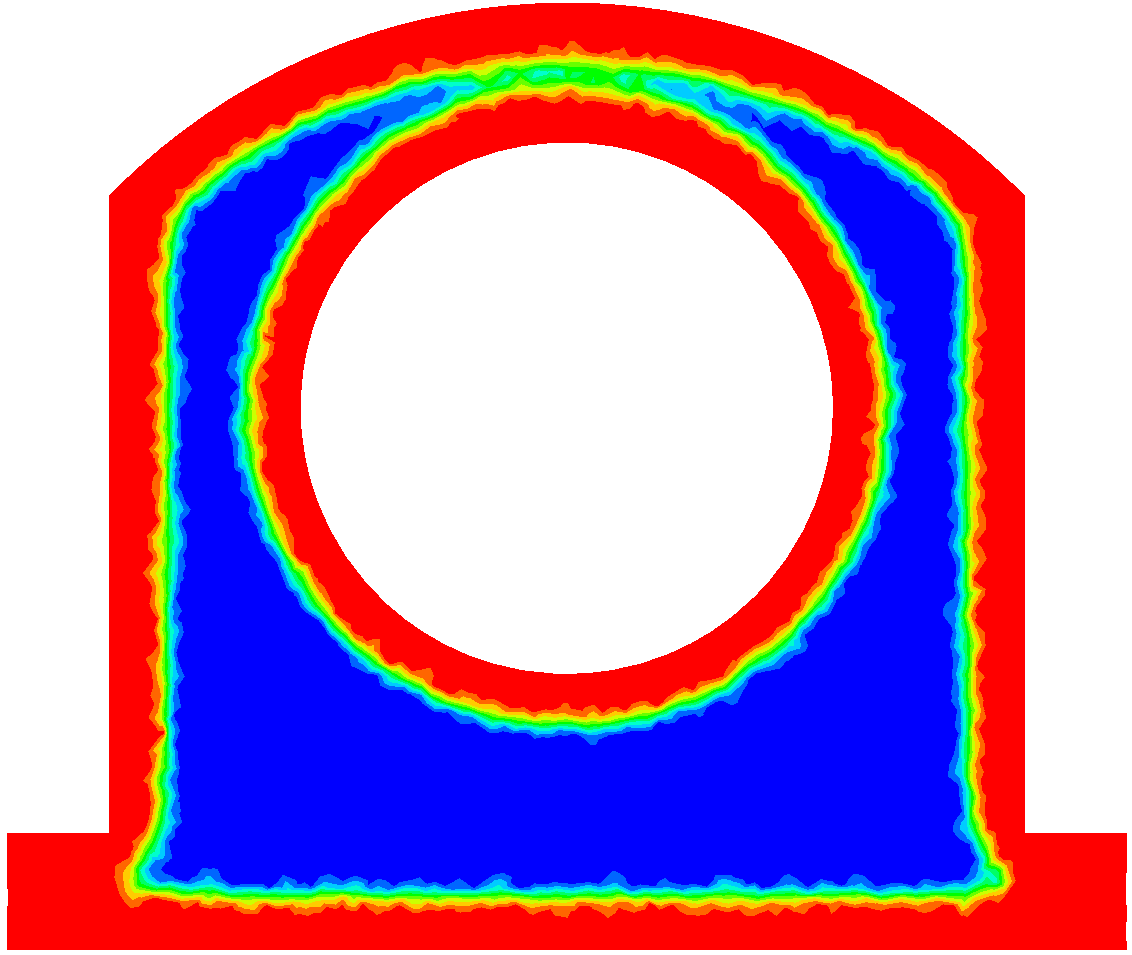}
		\caption{Time: 0.480 s}
	\end{subfigure}
	\begin{subfigure}{0.23\textwidth}
		\includegraphics[width=\textwidth]{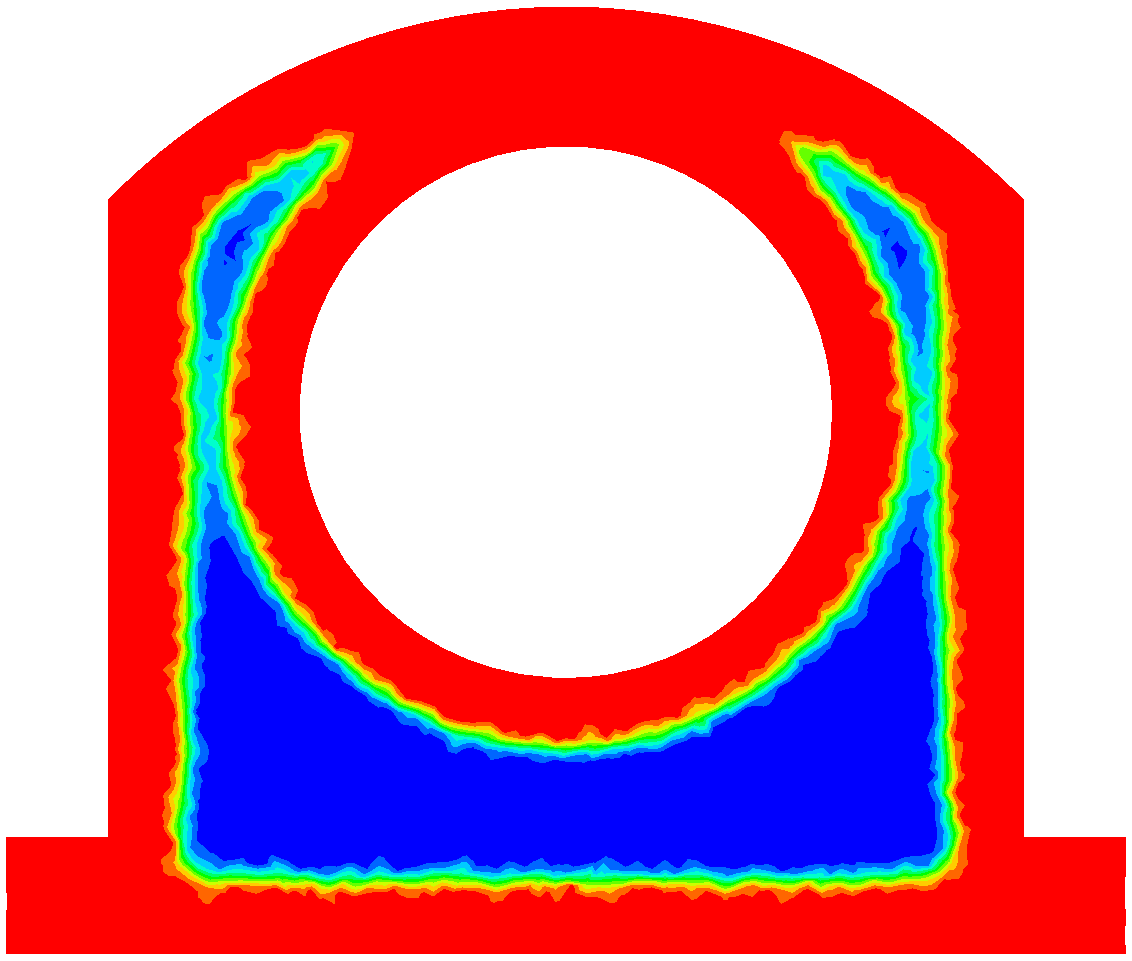}
		\caption{Time: 0.848 s}
	\end{subfigure}
	\begin{subfigure}{0.23\textwidth}
		\includegraphics[width=\textwidth]{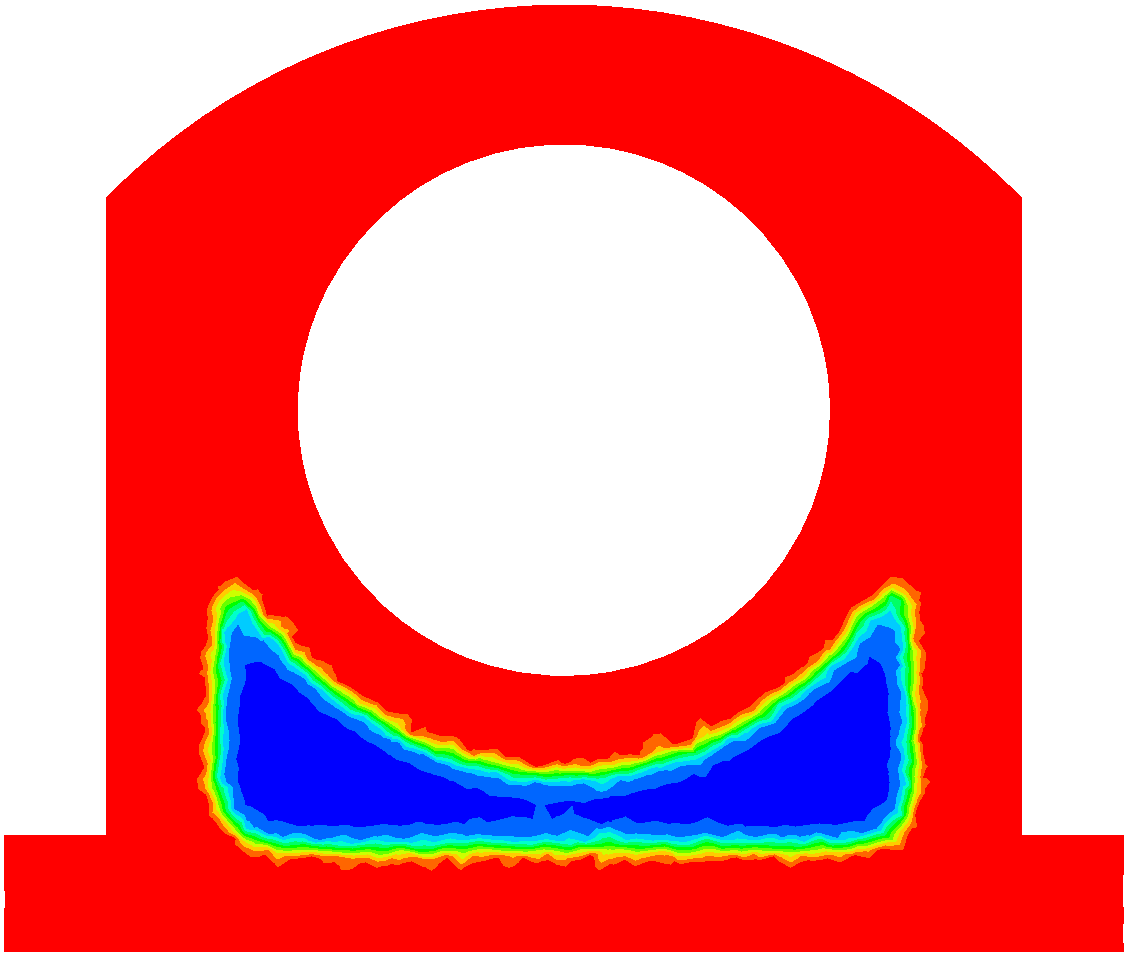}
		\caption{Time: 1.42 s}
	\end{subfigure}	
	\begin{subfigure}{0.05\textwidth}
		\includegraphics[width=\textwidth]{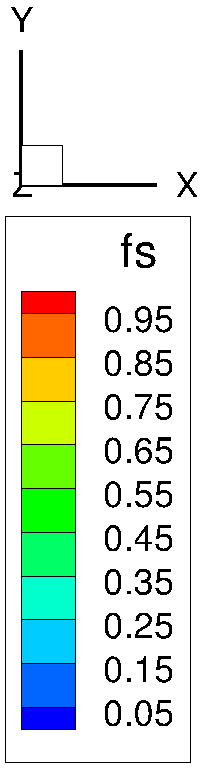}
	\end{subfigure}	
	\caption{Clamp Solid Fraction Contours (Solidification Time: 2.40 s)}
	\label{Fig:clamp_fs}
\end{figure}
\begin{figure}[H]
	\begin{subfigure}{0.325\textwidth}
		\includegraphics[width=\textwidth]{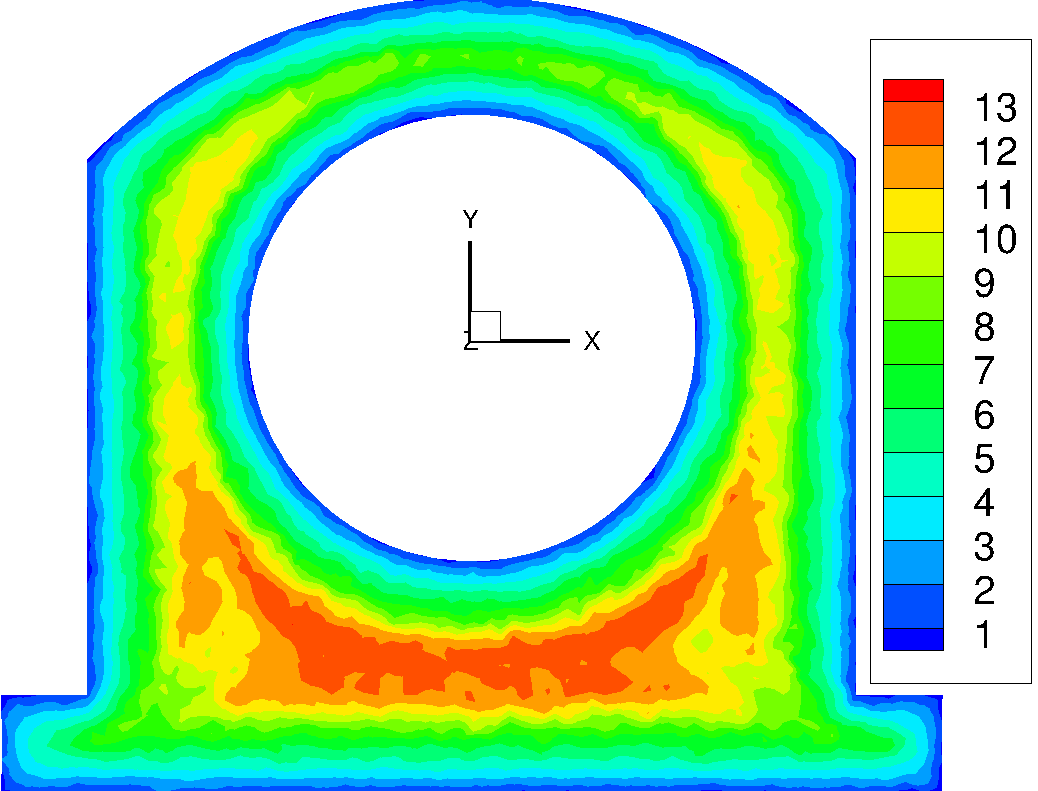}
		\caption{SDAS ($\mu$m)}
		\label{Fig:clamp SDAS}
	\end{subfigure}	
	\begin{subfigure}{0.325\textwidth}
		\includegraphics[width=\textwidth]{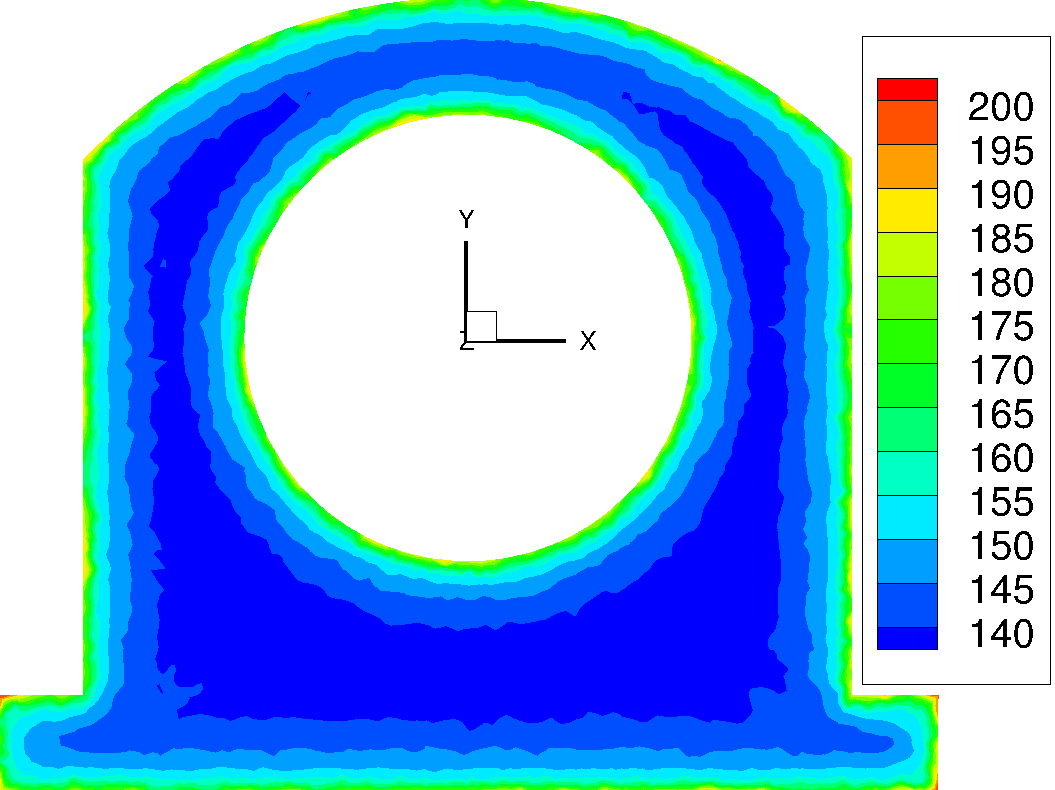}
		\caption{Yield Strength (MPa)}
		\label{Fig:clamp yield}
	\end{subfigure}
	\begin{subfigure}{0.325\textwidth}
		\includegraphics[width=\textwidth]{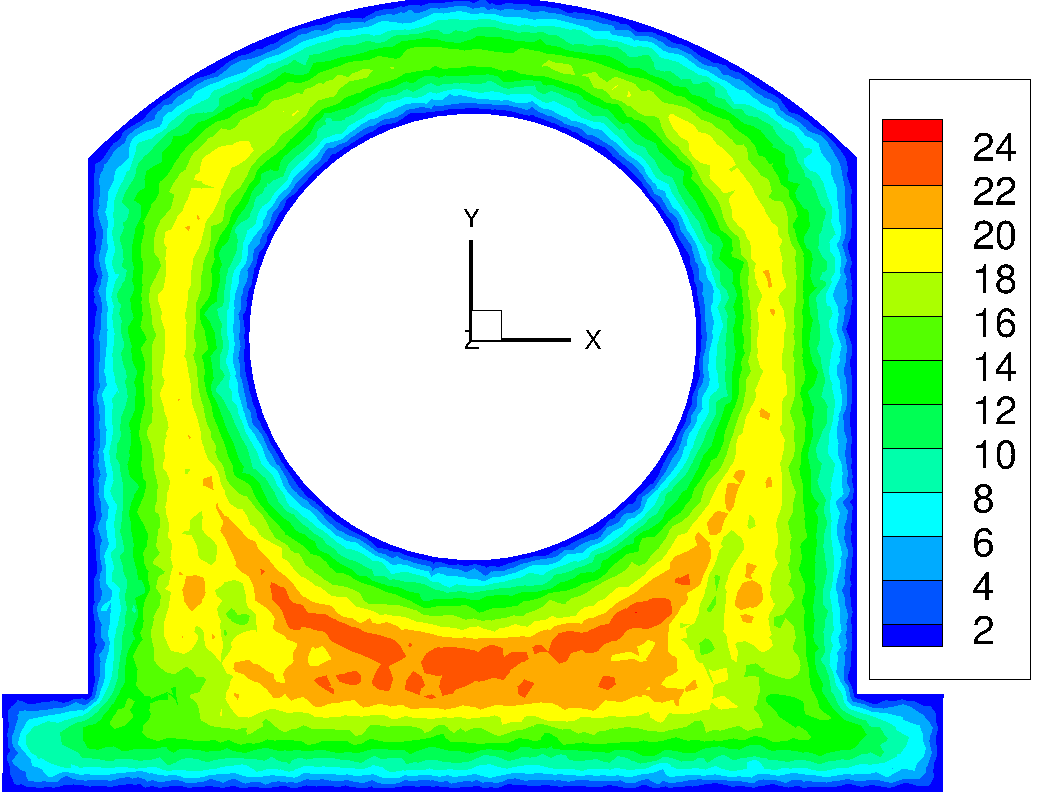}
		\caption{Grain Size ($\mu$m)}
		\label{Fig:clamp grain size}
	\end{subfigure}
	\caption{Clamp Microstructure Parameters}
	\label{Fig:clamp_microstruc}
\end{figure}
\par
\begin{figure}[H]
	\begin{subfigure}{0.325\textwidth}
		\includegraphics[width=\textwidth]{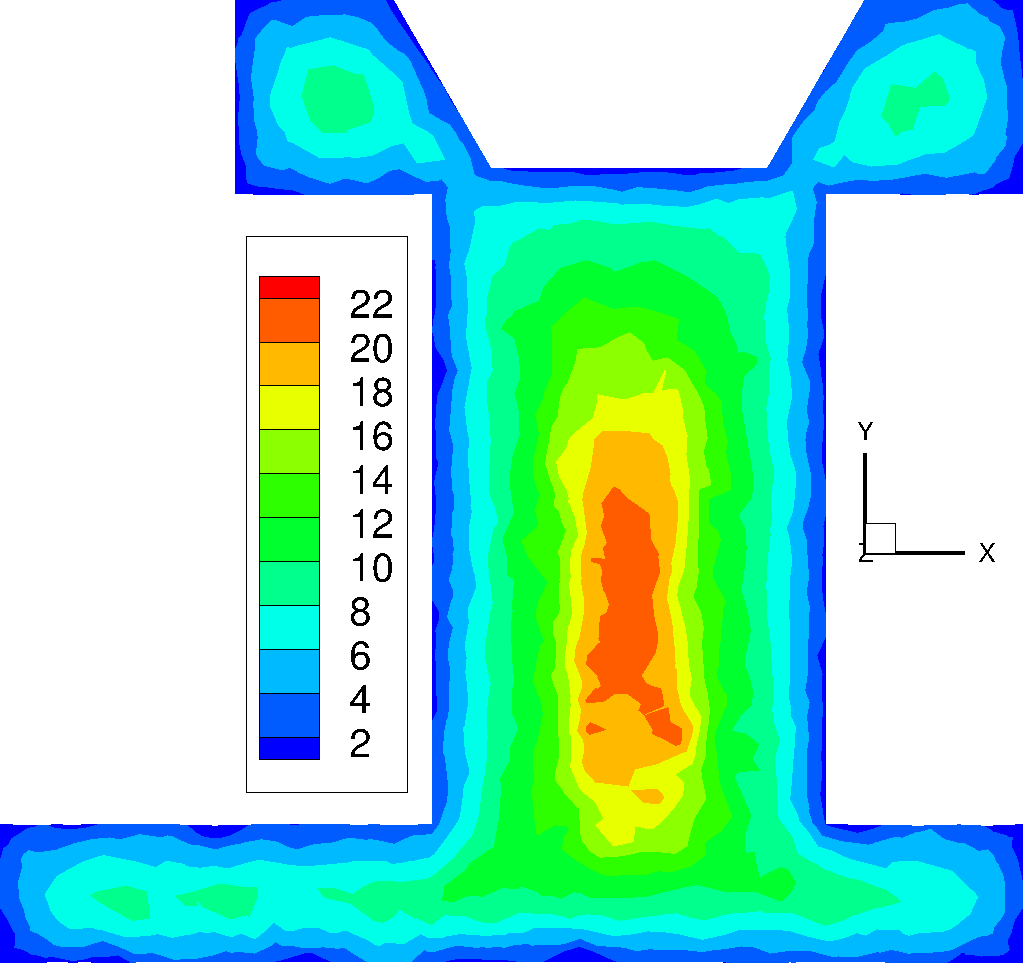}
		\caption{SDAS ($\mu$m)}
		\label{Fig:pulley SDAS}
	\end{subfigure}	
	\begin{subfigure}{0.325\textwidth}
		\includegraphics[width=\textwidth]{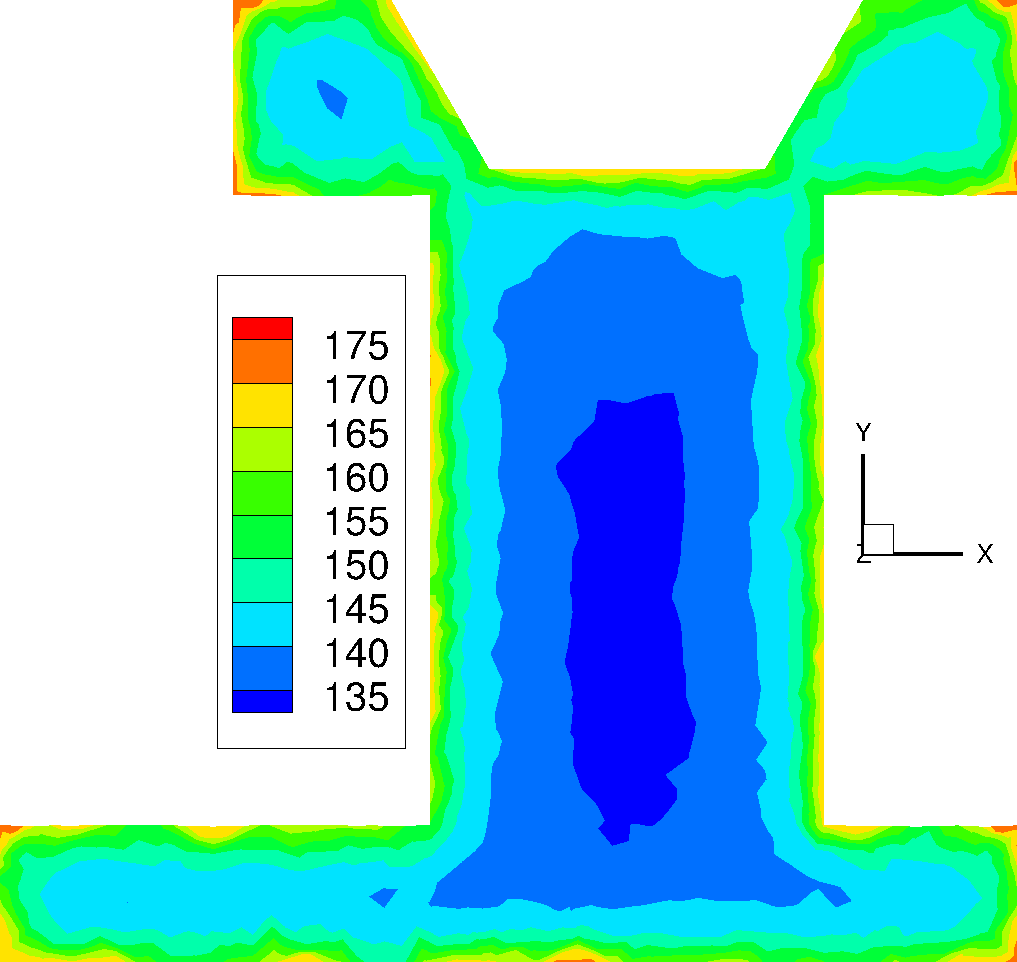}
		\caption{Yield Strength (MPa)}
		\label{Fig:pulley yield}
	\end{subfigure}
	\begin{subfigure}{0.325\textwidth}
		\includegraphics[width=\textwidth]{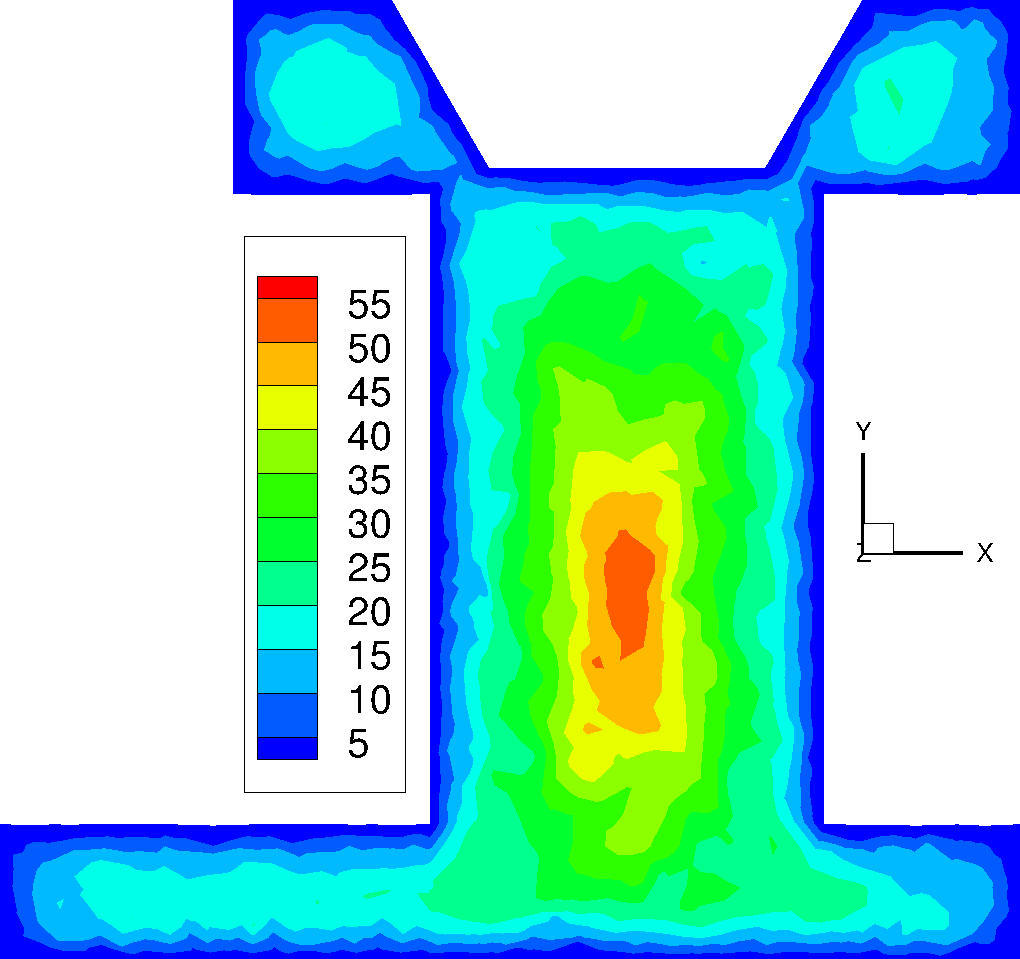}
		\caption{Grain Size ($\mu$m)}
		\label{Fig:pulley grain size}
	\end{subfigure}
	\caption{Pulley Microstructure Parameters}
	\label{Fig:pulley_microstruc}
\end{figure}

 Figure~\ref{Fig:clamp_T} plots the temperature contours for the clamp geometry along the mid-plane in Z-direction at different time-steps during the solidification. Figure~\ref{Fig:clamp_fs} shows the corresponding solid fraction contours. Figure~\ref{Fig:clamp_microstruc} plots micro-structure parameters SDAS, grain size and yield strength estimated using the previously described empirical models. As time progresses, the cooling rates and temperature gradients decrease. Hence, the regions with highest thickness or the core take longer time to solidify. The core regions typically have higher SDAS, grain size and lower yield strength as the grains have longer time to grow. Such trends are seen in Fig.~\ref{Fig:clamp_microstruc}. The pulley is axisymmetric in geometry and loading conditions. Thus, SDAS, grain size and yield strength are plotted only along the XY plane as shown in Fig.~\ref{Fig:pulley_microstruc}. Even for this case, the core region has higher SDAS and grain size and lower yield strength.
\section{Uncertainty Quantification Results}
The clamp geometry is subsequently used to study the effect of stochastic variation in the process parameters on the outputs. The process parameters which can be tuned in a die casting foundry are alloy composition (the solute percentage), initial melt temperature and the wall temperature. It is assumed that these three parameters follow normal distribution ($\mathcal{N}$) with the following means ($\mu$)  and standard deviations ($\sigma$):
\begin{enumerate}
	\item Solute concentration $C_0 \sim \mathcal{N} (\mu = 10, \sigma = 0.2) $ wt \%
	\item Initial molten alloy temperature $T_{init} \sim \mathcal{N} (\mu = 1000, \sigma = 1\% \mu = 10)$ K
	\item Wall temperature $T_{wall} \sim \mathcal{N} (\mu = 500, \sigma = 1\% \mu = 5)$ K
\end{enumerate}
\begin{figure}[H]
	\centering
	\includegraphics[width=0.45\textwidth]{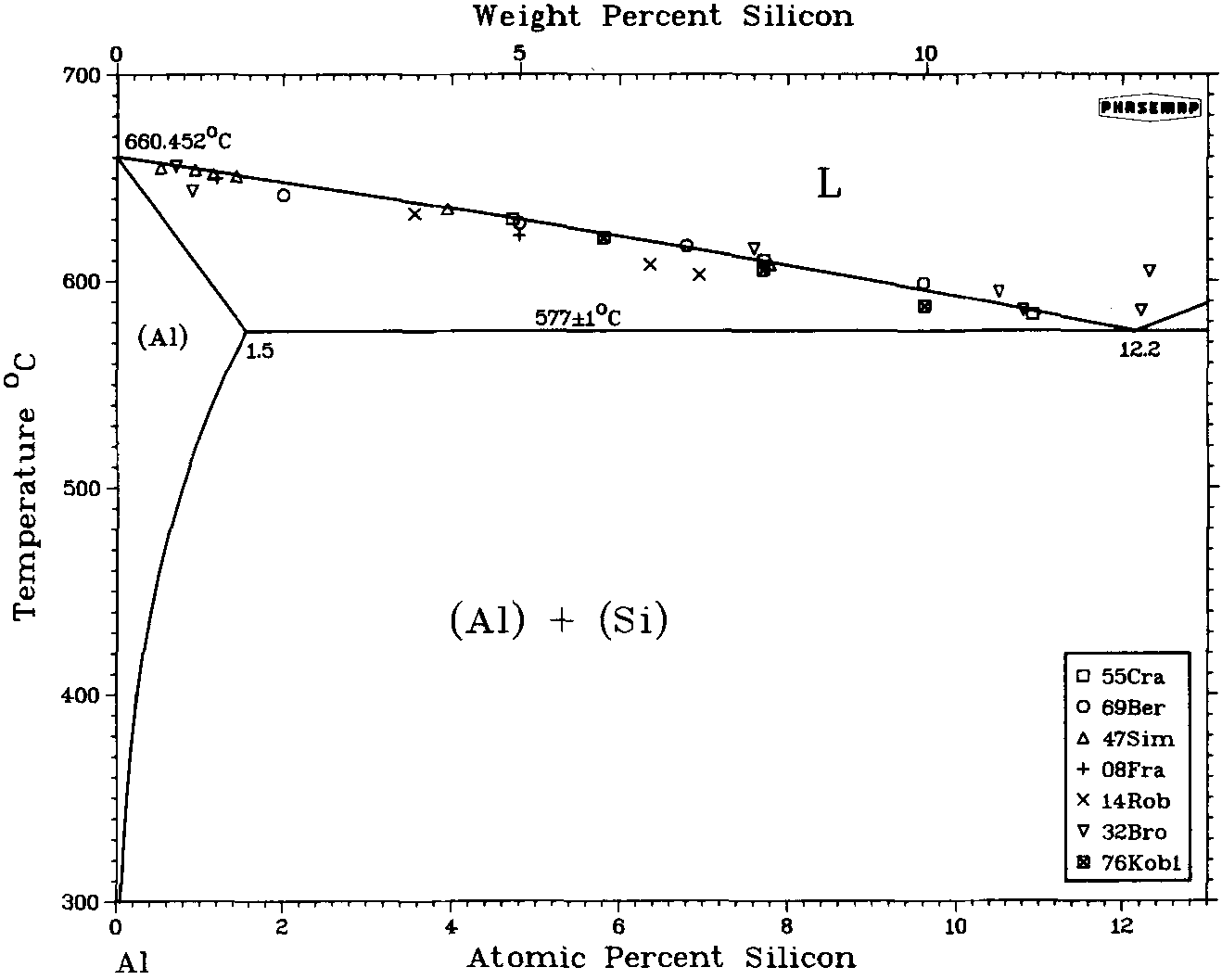}
	\caption{Al-Si Phase Diagram \cite{murray1984si}}
	\label{Fig:Al-Si Phase Diagram}
\end{figure}
Thermal properties like thermal conductivity, density and specific heat are estimated as a weighted linear combination of the properties of aluminum and silicon with solute concentration as the weight. Since there is a significant variation with temperature \cite{kaye1921tables}, temperature varying properties are used. Alloy properties are obtained from the Al-Si phase diagram (Fig.~\ref{Fig:Al-Si Phase Diagram}). The freezing and solidus temperatures are 660$^0$C and 577$^0$C, respectively. The solidus and liquidus curves can be approximated by straight lines. The partition coefficient can be estimated as the ratio of liquidus to solidus line slope \cite{dantzig2001modeling}. The liquidus temperature is calculated from the liquidus line as a function of the solute concentration. 
\begin{figure}[H]
	\centering
	\begin{tcbraster}[raster columns=2, raster equal height, 
		raster column skip=10pt, raster row skip=10pt, raster every box/.style={blank}]
		\tcbincludegraphics{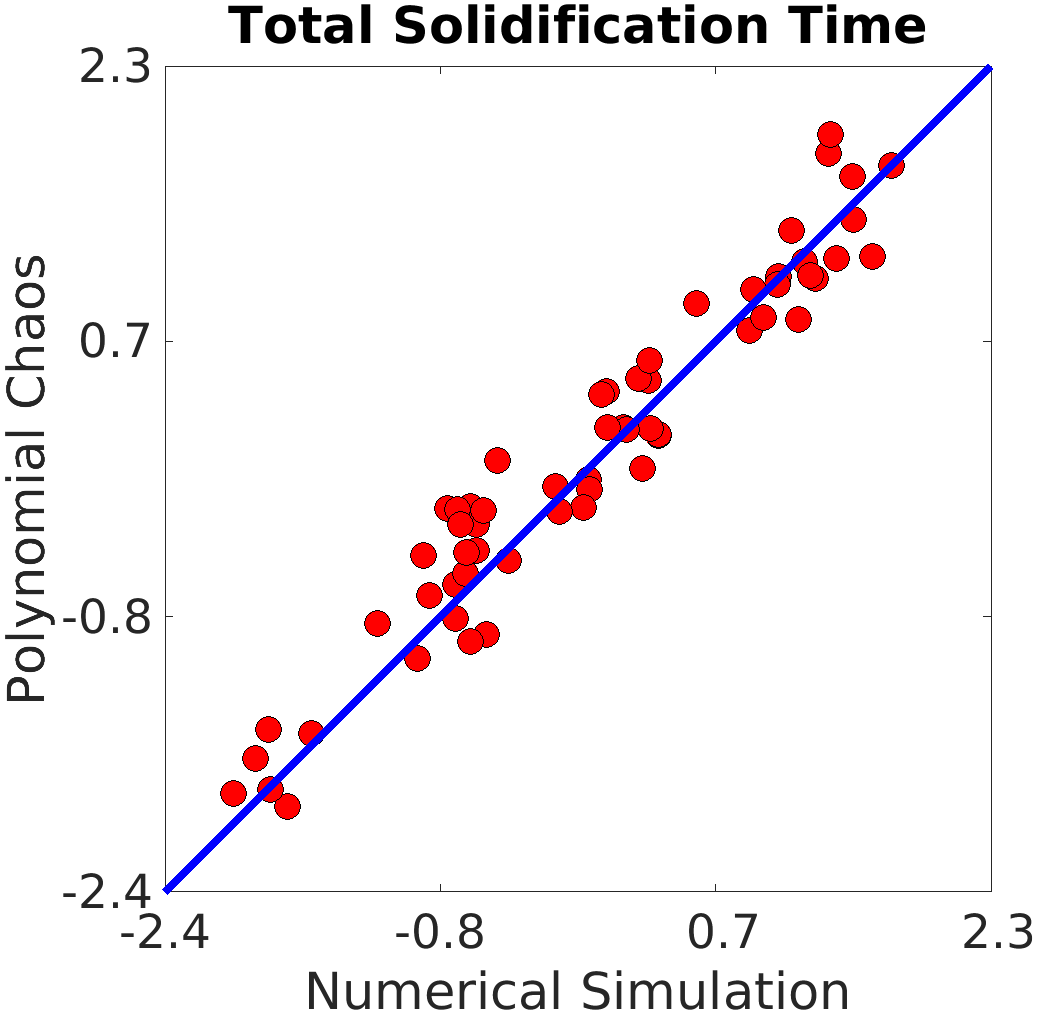}
		\tcbincludegraphics{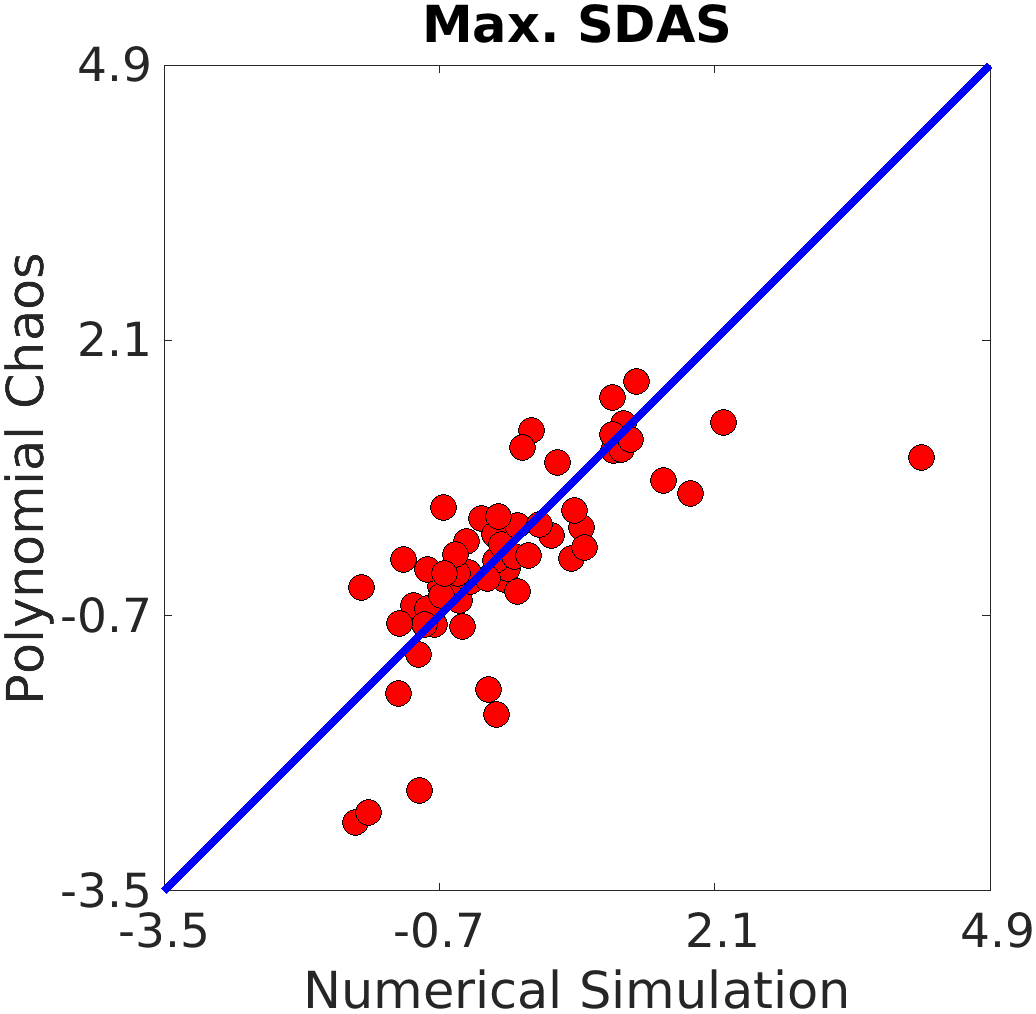}
		\tcbincludegraphics{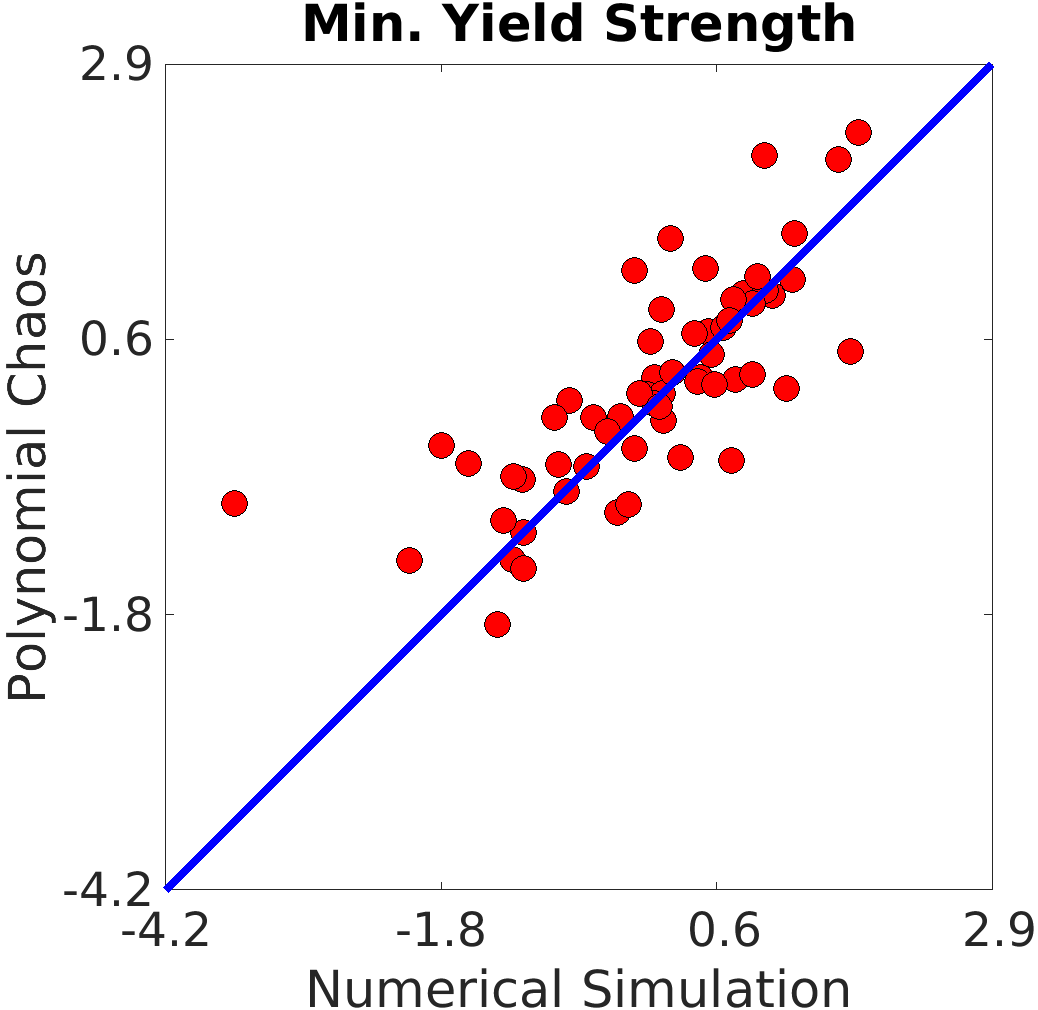}
		\tcbincludegraphics{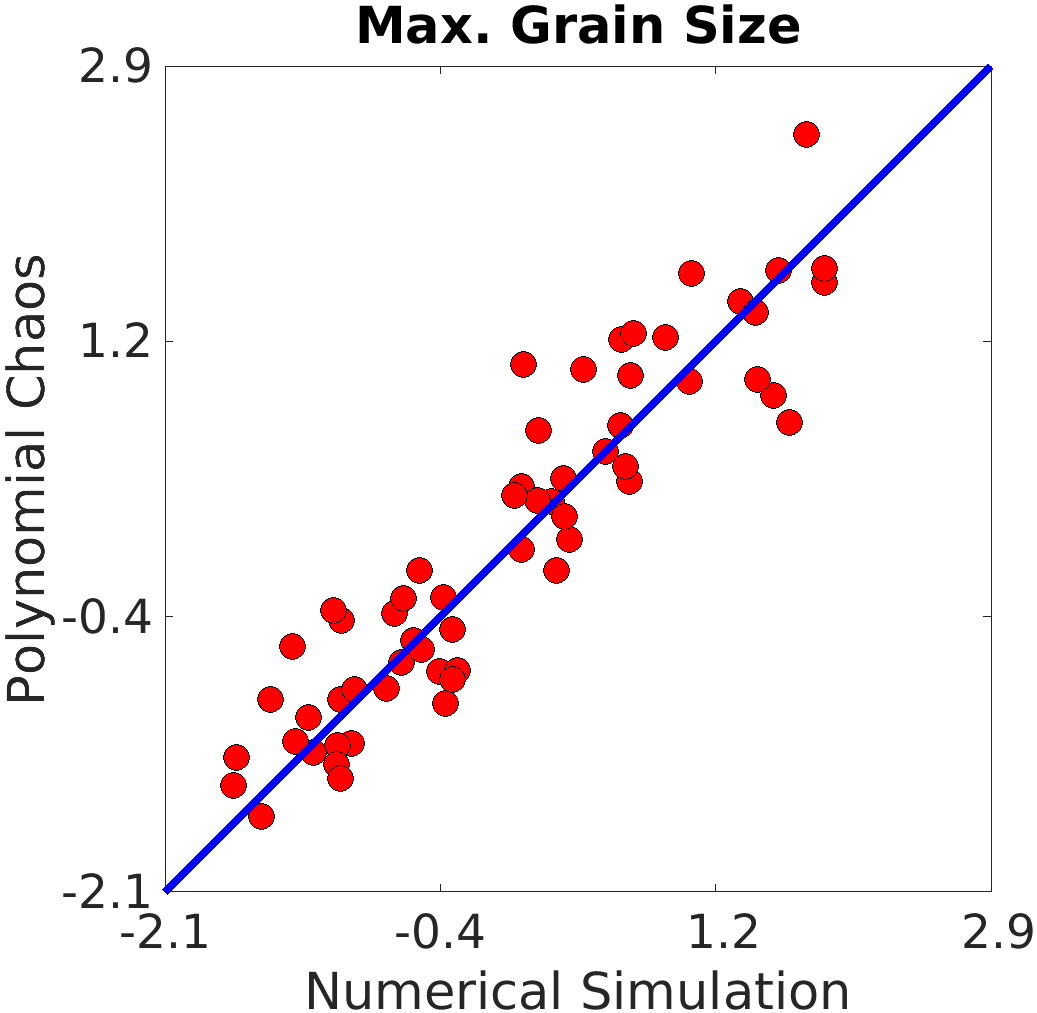}
	\end{tcbraster}
	\caption{Polynomial Chaos and Numerical Simulation Estimates (Non-dimensionalized)}
	\label{Fig:stoch coll error}
\end{figure}
\begin{table}[h]
	\begin{center}
		\begin{tabular}{|c|c|c|c|c|c|}
			\hline
			Accuracy Level&\# Sample Pts. &Output 1&Output 2&Output 3&Output 4\\
			\hline 
			3&27 & 1.20E-2 & 3.82E-2 & 2.70E-3 & 1.61E-2\\
			6&216 & 7.41E-3 & 3.35E-2 & 2.42E-3 & 1.42E-2\\
			\hline
		\end{tabular}
		\caption{Stochastic Collocation Convergence Analysis}
		\label{Table:Stochastic Collocation Convergence Analysis UQ}
	\end{center}
\end{table}
\par The impact of the three input parameters is studied on the following four outputs:
\begin{enumerate}
	\item Total solidification time
	\item Maximum SDAS
	\item Minimum yield strength
	\item Maximum grain size
\end{enumerate}
where, the maximum or minimum is taken over the entire domain. The polynomial chaos method with stochastic collocation described in section~\ref{Sec:Parameter Uncertainty Quantification} is used for uncertainty propagation. The error due to stochastic interpolation is estimated using 60 Latin hypercube samples. Estimates of the same output are obtained from the complete simulation and the polynomial chaos expansion (PCE) independently. The difference between these quantities non-dimensionalized by their mean value is defined as the error. Since the errors in estimating all the four outputs reduce by using higher accuracy levels (Table~\ref{Table:Stochastic Collocation Convergence Analysis UQ}), it can be seen that the collocation method has converged. The results presented correspond to the accuracy level 6. In order to get an idea of the interpolation error visually, both the simulation and PCE estimates can be plotted on the same graph. In the hypothetical scenario of exact interpolation, all the points should lie on the $Y=X$ line. However, there is a deviation from the line because of the interpolation error. Figure~\ref{Fig:stoch coll error} shows plots for each of the four outputs. Each output is normalized by subtracting its mean and dividing by its standard deviation and thus, is non-dimensional. Most of the points follow the trend of the $Y=X$ line except some outliers. The outliers correspond to those random samples which are too far from the means of the input parameters. Table~\ref{Table:Stochastic Collocation Convergence Analysis UQ} and Fig.~\ref{Fig:stoch coll error} prove that the polynomial chaos has converged and is accurate enough for further use.
\par
\begin{figure}[H]
	\centering
	\includegraphics[width=0.5\textwidth]{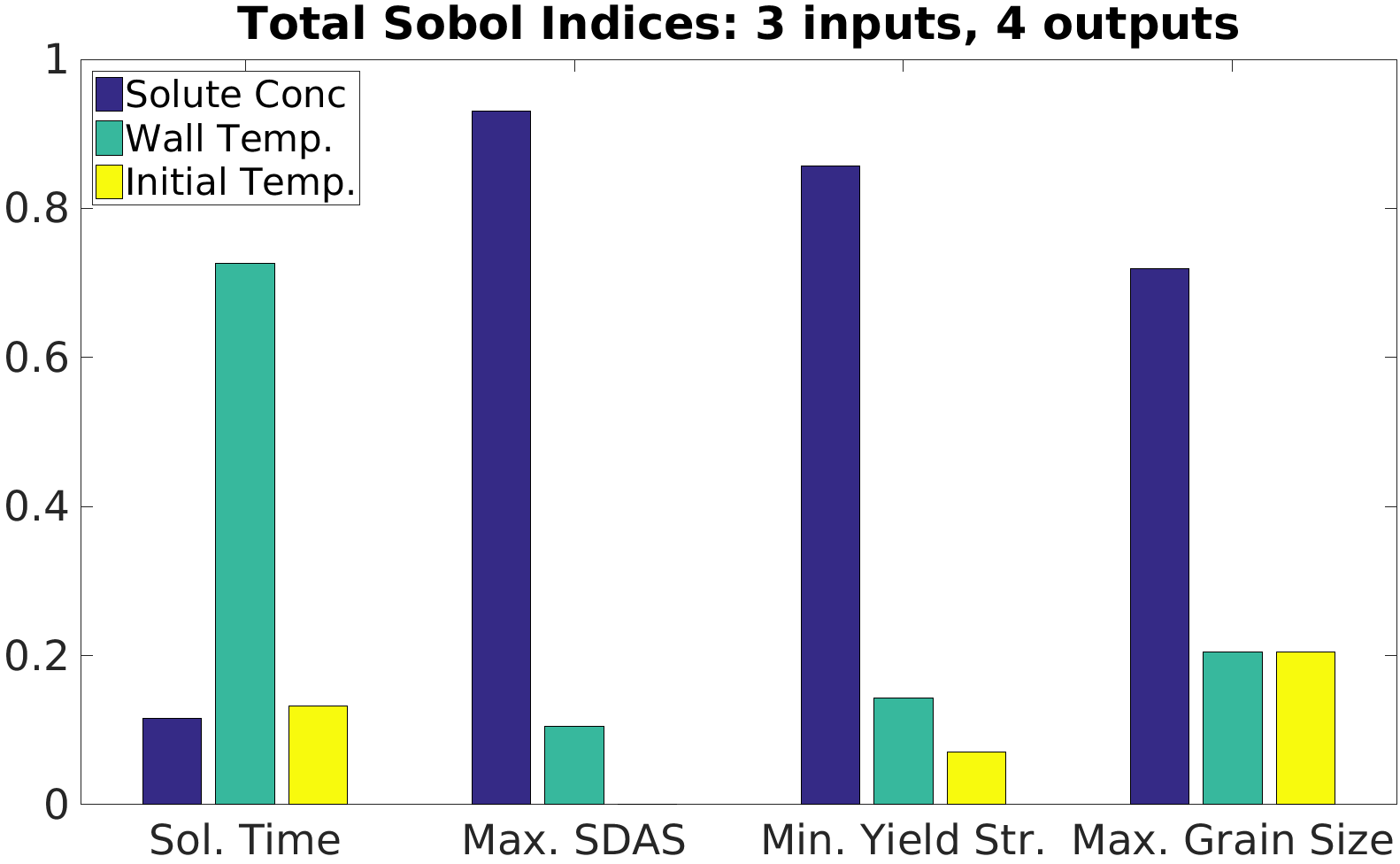}
	\caption{Sensitivity of 4 Outputs to 3 Inputs}
	\label{Fig:Sensitivity}
\end{figure}

 Sensitivity of each output with respect to each input can be easily estimated once an accurate polynomial chaos expansion is obtained. The sensitivity analysis tool of the software UQLAB \cite{marelli2014uqlab} is used to estimate the total Sobol indices (Fig.~\ref{Fig:Sensitivity}). It can be seen that each output is highly sensitive only to one input parameter. Sensitivity is practically important as it gives an idea as to which input parameter should be tightly controlled. Thus, other parameters can be loosely controlled saving cost but yielding desired product quality at the same time. The amount of heat extracted from the wall is proportional to the temperature gradient near the wall and hence, solidification time is highly sensitive to the wall temperature. On the other hand, the microstructure parameters like SDAS, yield strength and grain size depend on the solute concentration. Thus, these three parameters are more sensitive to solute concentration than the wall and initial temperatures.
\par
\begin{figure}[H]
	\centering
	\begin{tcbraster}[raster columns=2, raster equal height, 
		raster column skip=10pt, raster row skip=10pt, raster every box/.style={blank}]
		\tcbincludegraphics{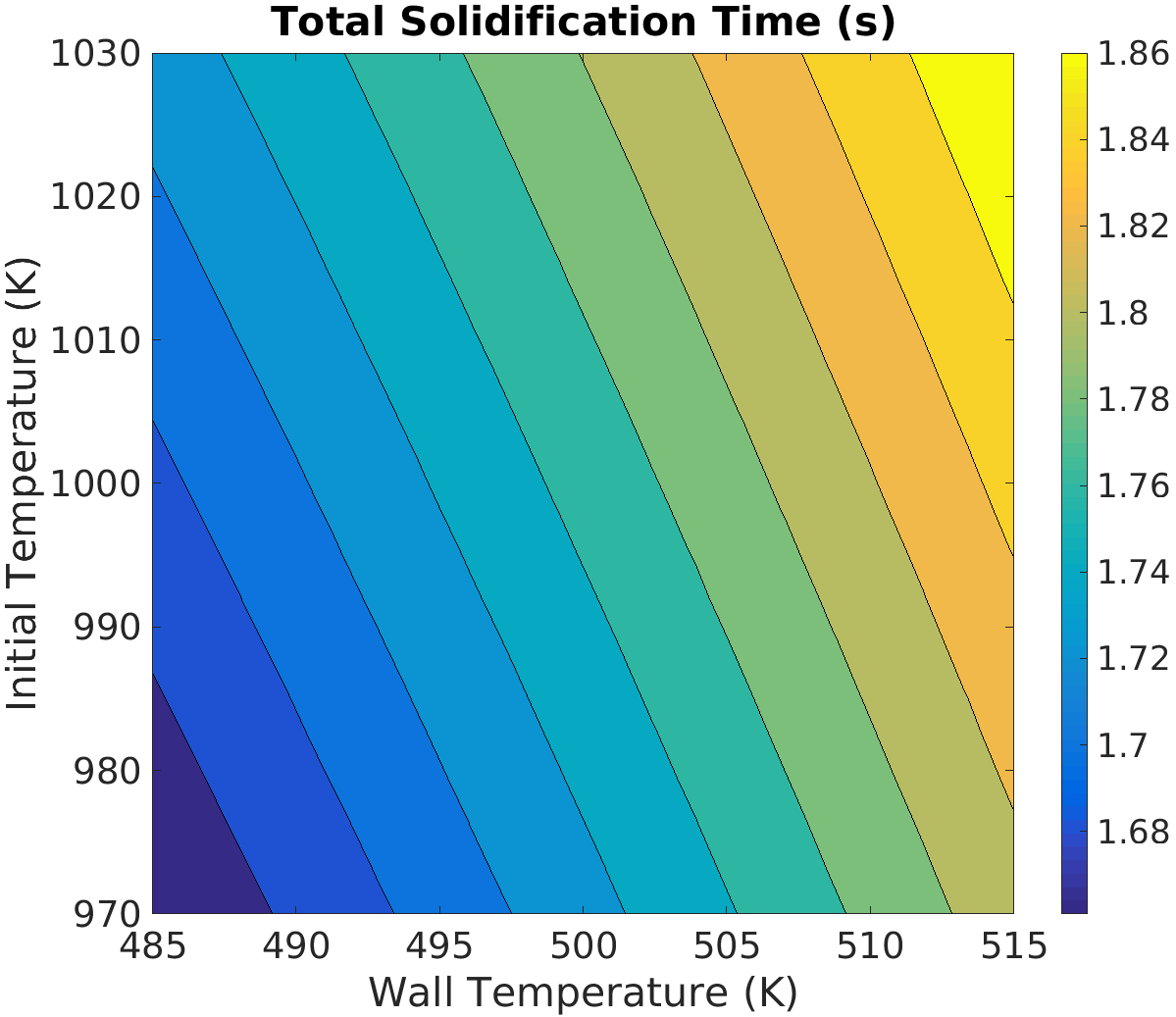}
		\tcbincludegraphics{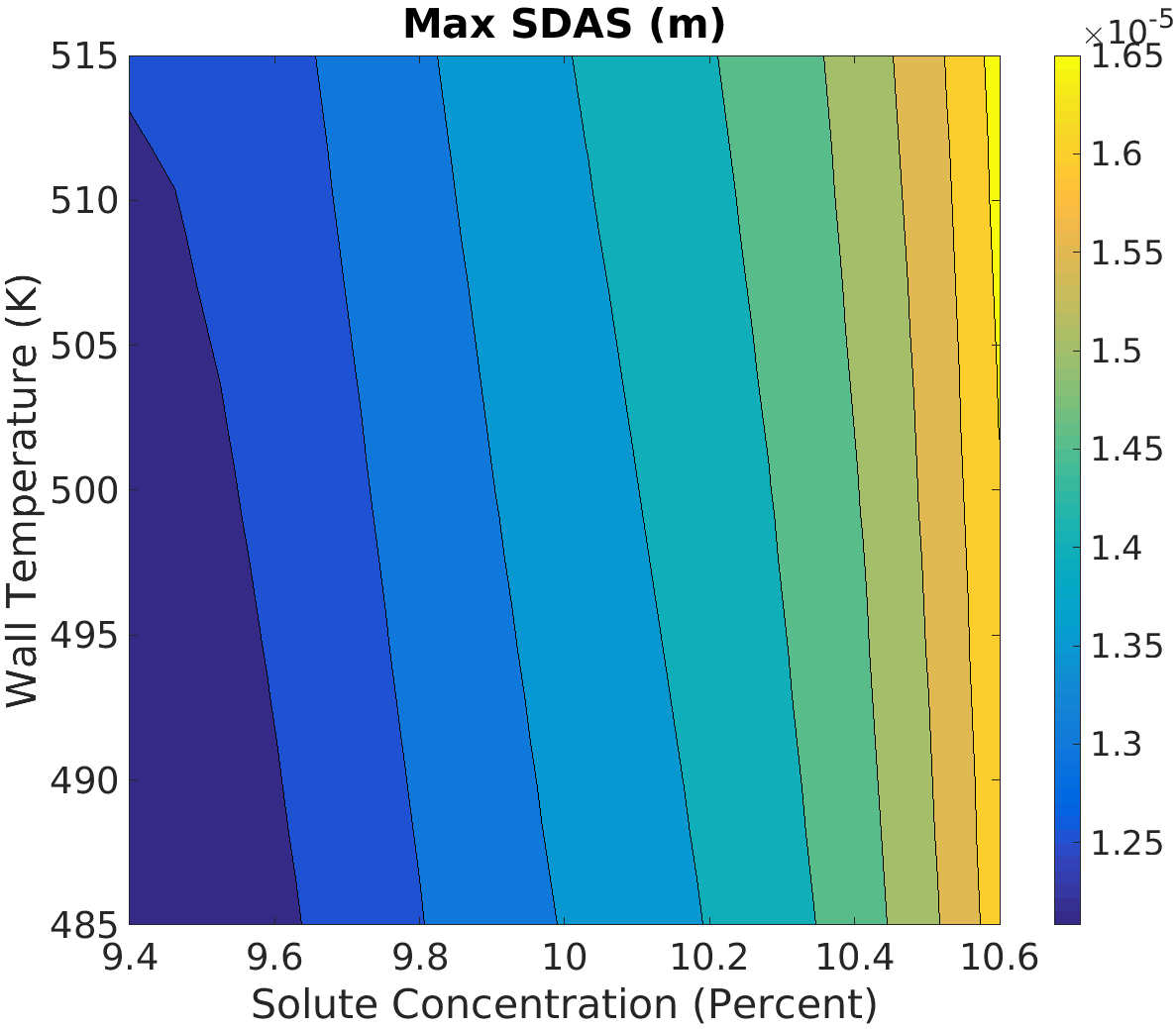}
		\tcbincludegraphics{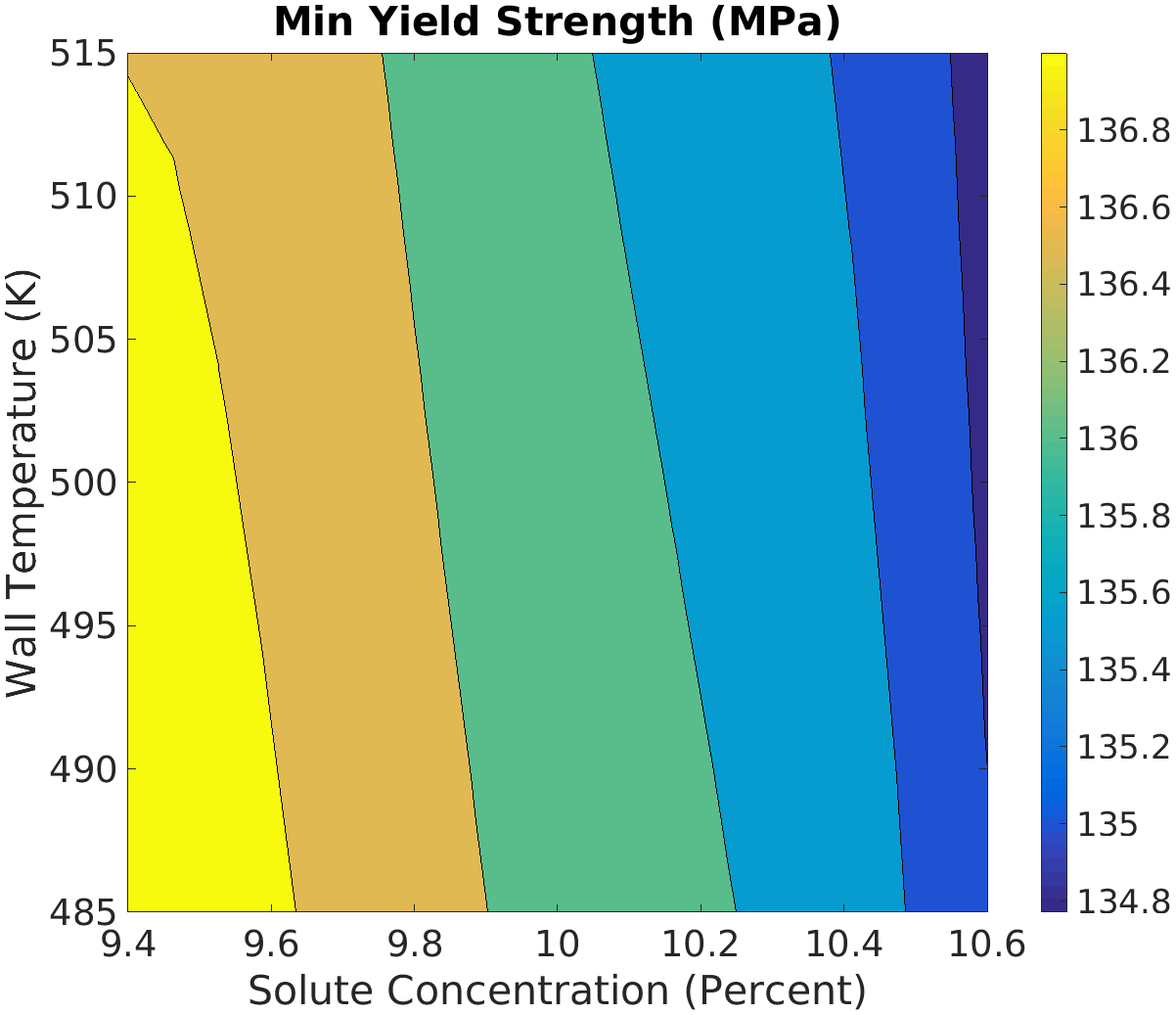}
		\tcbincludegraphics{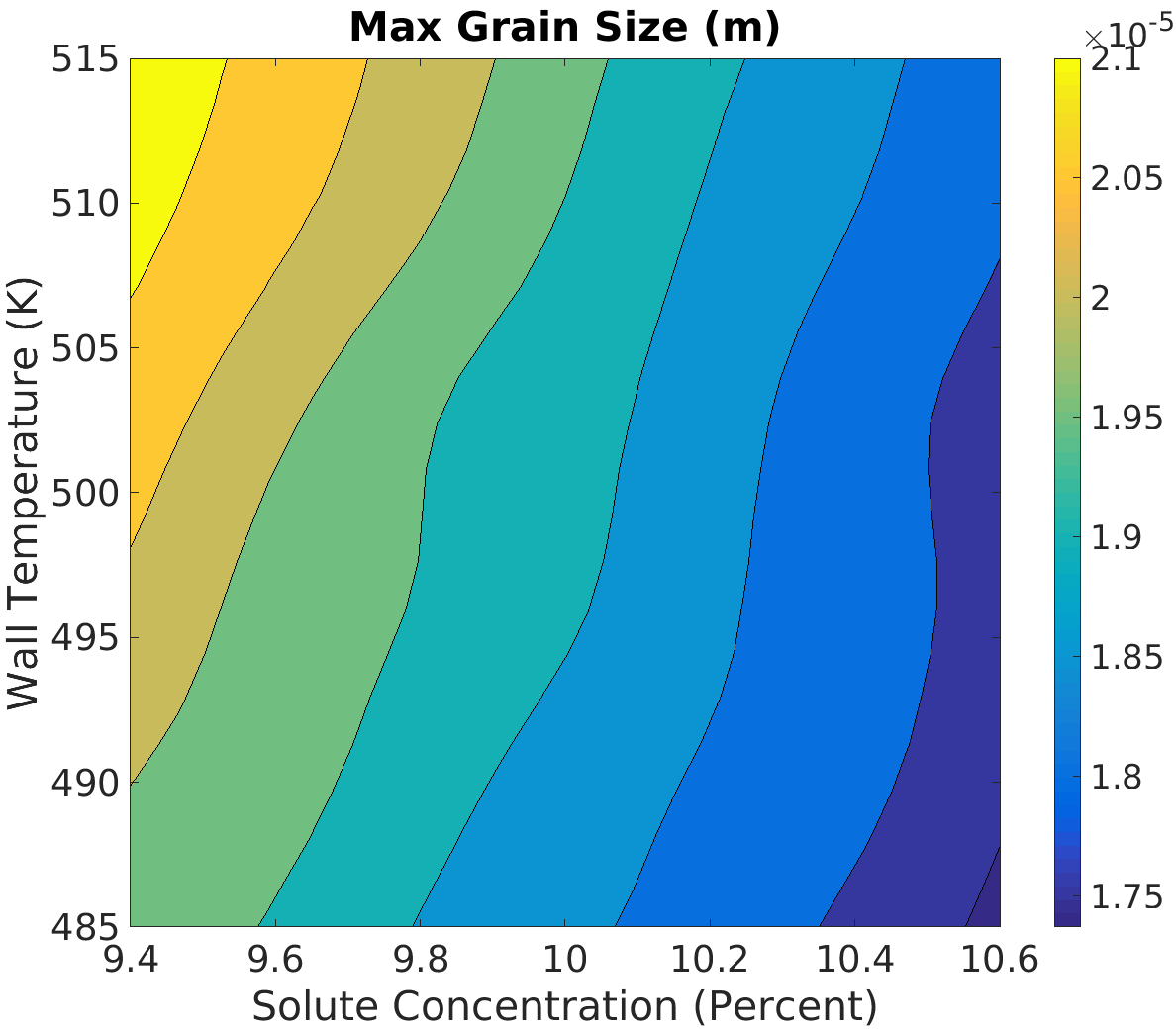}
	\end{tcbraster}
	\caption{Response Surfaces}
	\label{Fig:Response Surfaces}
\end{figure}

 In the response surface (Fig.~\ref{Fig:Response Surfaces}), the most important input is plotted on the x-axis. For solidification time response surface, the wall temperature (x-axis) and initial temperature (y-axis) are chosen. For other three outputs, solute concentration (x-axis) and wall temperature (y-axis) are chosen. It can be seen that all the response surface contours are nearly vertical. This confirms that the output is most sensitive to the input parameter plotted on the x-axis. The maximum grain size contours are non-linear. This implies that the sensitivity varies locally in the input parameter space. For other three outputs, the contours are almost linear and thus, it can be concluded that the local sensitivity is similar everywhere and independent of the input parameter value.  
\section{Conclusions}
This paper describes a numerical software framework OpenCast for simulations of solidification problems, including natural convection. Microstructure parameters and mechanical properties are estimated using published empirical relations. The flow equations are solved only in the liquid zone by modifying the coefficients of the discrete equations. The algebraic multigrid solver together with a Krylov subspace solver is used to solve the pressure Poisson equation. Complex geometries are meshed with unstructured hexahedral elements. Parameter uncertainty quantification is used as a wrapper over the deterministic simulations in order to assess the effect of stochastic variations in the inputs on the outputs.
\par The software is validated against published experimental results of solidification. This validation study shows the significance of stochastic analysis since it is observed that validation without uncertainty is unsuccessful. The validated software is used to simulate two practical die casting geometries. Sensitivity and uncertainty analysis is also performed. Sensitivity analysis shows that the product quality given by grain size and yield strength is highly sensitive to the solute concentration whereas, the productivity given by the solidification time is sensitive to the mold wall temperature. These results are practically useful as it gives an idea about the important input process parameters. The response surfaces show the variation of the outputs with the important input parameters. They can be used to get quick estimates of the outputs without running full deterministic simulations and also to get local sensitivities.
\par Although this paper demonstrates the ability of OpenCast to simulate die casting problems, it is a general purpose software. OpenCast can also be used to simulate other manufacturing processes such as sand casting, additive manufacturing, welding etc. The use of unstructured elements together with algebraic multigrid method adds complete flexibility to simulate arbitrary geometries. The coupling of uncertainty and sensitivity analysis tools enhance the power of the deterministic numerical method. 
\section*{Acknowledgment}
Authors acknowledge the financial support from Digital Manufacturing and Design Innovation Institute (DMDII Grant No. 15-07-06). The authors also thank Steve Udvardy and Beau Glim of North American Die Casting Association (NADCA) for providing industry contacts. Technical discussions with Alex Monroe of Mercury Castings and his suggestions were helpful in this work.
\section*{References}
\bibliography{References}
\end{document}